\documentclass[a4paper,12pt]{article}

\usepackage{amsfonts,amsthm,amsmath,amsfonts,latexsym,amssymb}
\usepackage[spanish, english]{babel}
\usepackage[ansinew]{inputenc}
\usepackage{graphicx,color,upref,float}
\usepackage[colorlinks=true,linkcolor=azul,citecolor=azul]{hyperref}

\newtheorem{fed}{Definition}[section]
\theoremstyle{definition}
\newtheorem{teo}[fed]{Theorem}
\newtheorem*{teo*}{Theorem}
\newtheorem{lem}[fed]{Lemma}
\newtheorem{cor}[fed]{Corollary}
\newtheorem{pro}[fed]{Proposition}

\theoremstyle{definition}

\newtheorem{rem}[fed]{Remark}

\newtheorem{exas}[fed]{Examples}

\newtheorem{defi}[fed]{Definition} 
\newtheorem{thm}[fed]{Theorem}

\def\la{\lambda}
\def\La{\Lambda}

\def\w{\omega}
\def\W{\Omega}

\definecolor{mylred}{rgb}{0.85,0.24,0.2}
\definecolor{myblue}{rgb}{0,0.33,0.55}
\definecolor{myyellow}{rgb}{0.42,0.24,0.52}
\definecolor{mygreen}{rgb}{0.12,0.5,0.29}
\definecolor{myred}{rgb}{0.74,0.13,0.13}
\definecolor{mylgreen}{rgb}{0.68,0.98,0.6}
\definecolor{mylyellow}{rgb}{0.86,0.85,0.55}
\definecolor{myllyellow}{rgb}{0.87,0.86,0.56}
\definecolor{naranja}{RGB}{249,153,96}
\definecolor{azul}{rgb}{0.1,0.6,0.86}
\definecolor{ble}{rgb}{0,0.33,0.75}

\def\noi{\noindent}

\def\EOE{\hfill $\blacktriangle$}

\def\bdem{\begin{proof}}
\def\edem{\end{proof}}

\def\bm{\left(\begin{array}}
\def\em{\end{array}\right)}
\def\ben{\begin{enumerate}}
\def\een{\end{enumerate}}
\def\barr{\begin{array}}
\def\earr{\end{array}}

\def\la{\lambda}
\def\La{\Lambda}

\def\N{\mathbb{N}}
\def\Z{\mathbb{Z}}
\def\R{\mathbb{R}}

\def\Q{\mathbb{Q}}

\def\cZ{\mathcal{Z}}

\newcommand{\peso}[1]{ \quad \text{ \rm  #1 } \quad }

\newcommand{\pint}[1]{\left \langle\, #1 \, \right\rangle}

\newcommand\myeq{\mathrel{\stackrel{\makebox[0pt]{\mbox{\normalfont\tiny (i)}}}{=}}}

\newcommand\myeqd{\mathrel{\stackrel{\makebox[0pt]{\mbox{\normalfont\tiny (ii)}}}{=}}}

\definecolor{mcolor}{rgb}{0,0.5,0.5}

\providecommand{\keywords}[1]{\textbf{\textbf{Keywords:}} #1}
\providecommand{\class}[1]{\textbf{\textbf{ AMS 2010 Mathematics Subject Classification:}} #1}

\begin{document}

\title{{\textbf{Tiling functions and Gabor orthonormal basis}}}

\author{Elona  Agora, Jorge Antezana, and Mihail N. Kolountzakis}

\date{}

\maketitle

\let\thefootnote\relax\footnotetext{The first author is partially supported by Grants MTM2013-40985-P, MTM2016-75196-P, and UBACyT 20020130100422BA. 
The second author is partially supported  by Grants CONICET-PIP 152, MTM2016-75196-P, and UNLP-11X585. }

\begin{abstract}

\noindent We study the existence of Gabor orthonormal bases with window the characteristic function of the set 
$\W=[0,\alpha] \cup [\beta+\alpha, \beta+1]$ of measure 1, with $\alpha, \beta>0$.  
By the symmetries of the problem, we can restrict our attention to the case $\alpha\leq1/2$. 
We prove that either if $\alpha <1/2$ or ($\alpha=1/2$ and $\beta\geq 1/2$) there exist such  Gabor orthonormal bases,  with window the characteristic function of the set $\W$, if and only if $\W$ tiles the line. Furthermore, in both cases, we completely describe the structure of the set of time-frequency shifts associated to these bases.

 \medskip

 \noindent \keywords{Spectral sets; Fuglede's Conjecture; Tilings; Packings; Gabor basis }  \\
\noindent \class{ Primary: 42C99, 52C22 }

\end{abstract}
 
\tableofcontents

\section{Introduction} 

In this note we study the Gabor orthonormal basis having as window a characteristic function of a union of two intervals on $\R$, of the form
\begin{equation} \label{W}
\W=[0,\alpha] \cup [\beta+\alpha, \beta+1],  
\end{equation}
and $|\W|=1$. By the symmetries of the problem, we can restrict our attention to the case where $\alpha\in(0,1/2]$ and $\beta>0$. 

\medskip

This paper draws on the ideas and some results given  in \cite{GLW}, where the structure of Gabor bases
with the window being the unit cube has been studied in general space of dimension $d$, producing surprising results even from dimension $d=2$. Here we restrict ourselves to space dimension $d=1$ but we vary the window in the simplest possible way: we examine the situation when the window is the indicator function of a union of two intervals satisfying~\eqref{W}.

\medskip

First we characterize the Gabor orthonormal basis
for a general function $g\in L^{2}(\R^d)$ in terms of a tiling condition involving 
the short time Fourier transform (of the window with respect to itself)
also known as Gabor transform (see Theorem \ref{Gonb char}).
This condition is completely analogous to the characterization of domain spectra, in the study of the so-called
Fuglede conjecture (see, for instance, \cite[Theorem 3.1]{K01}).
Our main result is Theorem \ref{geom cond} which provides the geometric conditions that characterize  
the Gabor orthonormal basis of $\chi_{\W}$ where $\W$ is given by \eqref{W}, and $\alpha<1/2$ or ($\alpha=1/2$ and $\beta\geq 1/2$). 

\medskip

As a consequence of Theorem \ref{geom cond} we can relate the existence of Gabor orthonormal bases using $\chi_\W$ with tiling properties of the set $\W$. Recall that we say that a set $\W\subseteq \R^d$ tiles $\R^d$ if
$$
\sum_{\la\in\La} \chi_{\W}(x-\la)=1
$$
almost everywhere for some set $\La\subseteq \R^d$. In \cite{L}  $\L$aba proved that for the union of two intervals the Fuglede conjecture is true. In particular, the tiling condition implies the existence of an orthonormal basis $\{e^{2\pi i \la \cdot}\}_{\la\in\La}$ for $L^2(\W)$. The set of frequencies $\La$ is called a spectrum for $\W$. Moreover, she proved that $\W$ tiles the real line, and therefore it has a spectrum, if and only if $\alpha$ and $\beta$ satisfy one of the following two conditions:
\begin{itemize}
\item[(i)]  $\displaystyle 0<\alpha < 1/2 $ and  $\beta \in \N$;
\item[(ii)] $\alpha =1/2$  and $\displaystyle \beta \in \frac{1}{2}\N$.
\end{itemize} 
The sets $\W$ given by (i)  are  (up to translation and reflection) the $\Z$-tiles consisting of two intervals. On the other hand, the sets in the case (ii) are also $\Z$-tiles if $\beta$ is integer, otherwise, if $\displaystyle \beta=\frac{n-1}{2}$, for any even $n\in \N$ these sets are tiles with 
respect to the set of translations
$$
n\Z \cup \left(n\Z +\frac{1}{2}\right) \cup  \dots \cup \left(n\Z +\frac{n-1}{2}\right).
$$
By  $\L$aba's result, note that either in the case (i) or in the case (ii) there exists Gabor orthonormal basis
$$
\{e^{2\pi i \nu \cdot}\chi_\W(\cdot-t): (t,\nu)\in\Gamma\}
$$
where $\Gamma$ can be taken as a product of the set used to tile $\R$ with $\W$, and any spectrum $\La$ for $\W$. A natural question is whether or not there are other sets $\W$ as in \eqref{W} that produce Gabor orthonormal basis. As a consequence of Theorem \ref{geom cond} we get that, among the sets $\W$ that satisfy $\alpha<1/2$ or ($\alpha=1/2$ and $\beta\geq 1/2$),  only those which tile the real line can produce Gabor orthonormal bases (Theorem \ref{caso menor}). Moreover, in these cases we can describe the structure of the sets of time-frequency translations that produce those orthonormal bases (Theorems \ref{necstandard} and \ref{necstandard2}).

\medskip

One of the basic tools in this characterization is the behavior of the zero set of the short time Fourier transform of the characteristic function of $\W$,
a situation also observed in \cite{GLW} and which closely mirrors what happens in the related Fuglede problem \cite{K01}. Even though, this zero set is completely described in Appendix \ref{Zero set}, the aim of Appendix \ref{Zero set} is complementary and it is not directly used in the proof of the main result.

\section{Preliminaries}

Given $g\in L^2(\R^d)$, and $\la=(t,\nu)\in\R^{d}\times\R^d$, let $g_\la$ denote the \textbf{time-frequency shift} of $g$ defined by
$$
g_\la(x)= g(x-t)e^{2\pi i \pint{x, \nu}}.
$$
From now on, $\R^d\times\R^d$ will be identified with $\R^{2d}$. Let $\La$ be a discrete countable set of $\R^{2d}$. 
The \textbf{Gabor system} associated to the so-called window $g$ consists of the set of time-frequency shifts $\{g_{\la}\}_{\la\in \La}$. 
We are especially interested in the case when a Gabor system is an orthogonal basis (a Gabor basis).

\subsection{Short time Fourier transform} 

In this subsection we will recall the definition and some properties of the short time Fourier transform, also known as Gabor transform. 
More information can be found in \cite{Gro}.
\begin{fed} 
Given $g\in L^2(\R^d)$, the \textbf{short time Fourier transform} with \textbf{window} $g$ is defined for any $f\in L^2(\R^d)$ by
\begin{equation} \label{Vg}
V_gf(t,\nu):= \int_{\R^{d}} f(x) \overline{g(x-t)}e^{-2\pi i \pint{x,\nu}}dx =\pint{f, g_{(t,\nu)}}.
\end{equation}
\end{fed}

One of the most useful properties of this transform is the following relation.

\begin{pro}\label{Orthogonality properties of the STFT} Let $f_1,f_2,g_1,g_2\in L^2(\R^d)$. 
Then, $V_{g_j}f_j\in L^2(\R^{2d})$ for $j=1,2$, and
$$
\pint{V_{g_1}f_1,V_{g_2}f_2}=\pint{f_1,f_2}\overline{\pint{g_1,g_2}}.
$$
\end{pro}

\begin{cor}\label{isometria}
If $f,g\in L^2(\R^d)$, then
$$
\|V_g f\|^2_2=\|f\|^2_2\|g\|^2_2.
$$
In particular, if $\|g\|_2=1$ then $V_g$ is an isometry from $L^2(\R^d)$ into $L^2(\R^{2d})$.
\end{cor}

Another consequence is the following result.

\begin{cor}\label{denso}
Let $f,g\in L^2(\R^d)$ such that $g\neq 0$ in $L^2$. 
If $\pint{f,g_\w}=0$ for every $\w\in\R^{2d}$, then $f=0$ in $L^2$.
\end{cor}

\medskip

We conclude this section with the following proposition that provides the behavior of 
Gabor transforms  with respect to time-frequency shifts.

\medskip

\begin{pro}\label{TF shift y la Vg}
Given $f,g\in L^2(\R^d)$ and $\la=(t,\w)\in\R^{2d}$
$$
V_g \,f_\la\,(x,\nu)= e^{-2\pi i \pint{t,\nu-\w}} V_g f(x-t,\nu-w).
$$
\end{pro}

\subsection{Tilings}\label{sec:tilings}

Let $\delta_\la$ be the unit point mass sitting at the point $\la \in \R^n$. 
Recall that, given a function $h \geq 0$, the convolution $h*\delta_\la(x)$ in the distribution sense gives the translation $h(x-\la)$. 
Let $\delta_\La$ denote the measure 
$$
\delta_{\La}=\sum_{\la\in \La} \delta_\la,
$$
where  $\La$ is a discrete set of $\R^n$.  Then, if 
$$
\delta_\La * h(x)=\sum_{\la \in \La} h (x-\la)=1,  
$$ 
for almost every $x\in\R^n$,  we say that $h$ \textbf{tiles} $\R^{n}$ with translation set $\La\subseteq \R^n$.
If, for example, $h= \chi_\W$, where $\W$ is a measurable subset of $\R^n$ then the tiling condition means that 
$$
\la+\W:= \{\w+\la:\, \w\in \W \}, \peso{where} \la \in \La
$$ 
intersect each other in a set of zero Lebesgue measure, and their union covers the whole space 
except perhaps for a set of measure zero. We often denote the situation $h * \delta_\Lambda = 1$ (tiling) by
$$
h+\Lambda = \R^d.
$$
Similarly we say that $h+\Lambda$ is a \textbf{packing}, and we denote this by
$$
h+\Lambda \le \R^d,
$$
if $h*\delta_\Lambda \le 1$ almost everywhere. If $h=\chi_\W$, we also use the notation
$$
\W+\La \leq \R^d
$$
for packing and
$$
\W+\La = \R^d
$$
for tiling.

\section{Orthogonality and the zero set of Gabor transform}

Throughout this section $\La$ will denotes a discrete subset of $\R^{2d}$. Given a Borel set $\W$, which has finite measure, the orthogonality of exponential families in $L^2(\W)$ can be studied using the Fourier transform of $\chi_\W$. In this section we will show that there exists a similar connection between the orthogonality properties of a Gabor system $\{g_{\la}\}_{\La}$ and the zero set of $V_gg$.

\begin{thm} \label{Gonb char}
Given $g\in L^2(\R^d)$ and $\|g\|_2=1$, the following statements are equivalent:
\begin{itemize}
\item[(i)]  The Gabor system $\{g_{\la}\}_{\La}$ forms an orthonormal system.
\item[(ii)] 
$
\displaystyle \sum_{\la\in \La} |V_gf(\w-\la)|^2\leq1, 
$
for every $f\in L^2(\R^d)$ so that $\|f\|_2=1$, and $\w \in \R^{2d}$.
In other words
$$
|V_gf|^2 + \Lambda\ \ \ \mbox{is a packing}.
$$
\item[(iii)] 
$
\displaystyle \sum_{\la\in \La} |V_gg(\w-\la)|^2 \leq 1,
$ 
for every  $\w \in \R^{2d}$.
In other words
$$
|V_gg|^2+\Lambda\ \ \ \ \mbox{is a packing}.
$$
\end{itemize}
We also have that the following are equivalent:
\begin{itemize}
\item[(iv)]
The Gabor system $\{g_{\la}\}_{\La}$ forms an orthonormal basis.
\item[(v)]
$
\displaystyle \sum_{\la\in \La} |V_gf(\w-\la)|^2 = 1, 
$
for every $f\in L^2(\R^d)$ so that $\|f\|_2=1$, and $\w \in \R^{2d}$.
In other words
$$
|V_gf|^2 + \Lambda\ \ \ \mbox{is a tiling}.
$$
\item[(vi)] 
$
\displaystyle \sum_{\la\in \La} |V_gg(\w-\la)|^2 = 1,
$
for every  $\w \in \R^{2d}$.
In other words
$$
|V_gg|^2+\Lambda\ \ \ \ \mbox{is a tiling}.
$$
\end{itemize}
\end{thm}

\bdem

Let us first prove the equivalence of (i), (ii) and (iii). Clearly (ii) $\Rightarrow$ (iii). Hence, it is enough to show that (i)  $\Rightarrow $ (ii) and (iii)  $\Rightarrow $ (i).

\bigskip

(i)  $\Rightarrow $ (ii)  Assume that $\{g_{\la}\}_{\la\in \La}$ forms an orthonormal system for $L^2(\R^d)$. Since
\begin{equation}\label{eq symmetry}
|\pint{g_{\la_1},g_{\la_2}}|=|\pint{g_{-\la_1},g_{-\la_2}}|,
\end{equation} 
the Gabor system $\{g_{-\la}\}_{\la\in \La}$ is also orthonormal. 
Then, by the Bessel's inequality,
\begin{equation} \label{PI}
1\geq\sum_{\la\in \La} |\pint{f,g_{-\la}}|^2 =\sum_{\la\in \La} |V_gf(-\la)|^2,
\end{equation}
for every $f\in L^2(\R^d)$, so that $\|f\|_2=1$. By Proposition \ref{TF shift y la Vg}, if $\w\in\R^{2d}$ 
and replacing $f$ by $f_{-\w}$,  we get
\begin{equation*} 
1\geq\sum_{\la\in \La} |\pint{f_{-\w},g_{-\la}}|^2 =\sum_{\la\in \La} |V_gf(\w-\la)|^2. 
\end{equation*}

\medskip

(iii)  $\Rightarrow  $ (i) First of all, note that taking $w=\la_0$ for some $\la_0\in \La$ we obtain
$$
|V_g(0)|^2+\sum_{\la\neq \la_0} |V_gg(\la_0-\la)|^2 \leq 1.
$$
Since $V_gg(0)=\|g\|_2^2=1$, we get that for every $\la\in \La\setminus\{\la_0\}$
$$
0=|V_gg(\la_0-\la)|=|\pint{g_{\la_0},g_{\la}}|.
$$
Therefore, $\{g_{\la}\}_{\la\in\La}$ is an orthonormal system. 

\medskip

Now, assume (vi).
It follows from the equivalence of (i), (ii) and (iii) that the Gabor system
is orthonormal.
In order to prove that the Gabor system is complete it is enough to prove that
$$
\sum_{\la\in\La} |\pint{f,g_\la}|^2=1
$$
for a set of norm one elements $f\in L^2(\R^d)$ that is dense in the sphere of $L^2(\R^d)$. By \eqref{eq symmetry} and the hypothesis we know that for every $\w\in \R^{2d}$
$$
\sum_{\la\in\La} |\pint{g_{\w},g_{\la}}|^2=\sum_{\la\in\La} |\pint{g_{-\w},g_{-\la}}|^2=\sum_{\la\in \La} |V_gg(\w-\la)|^2=1.
$$
By Corollary \ref{denso} the family $\{g_\w\}_{\w\in\R^{2d}}$ is complete. So $\{g_{\la}\}_{\la\in\La}$ is an orthonormal basis. 
\edem

\medskip

As is the case of Fuglede problem \cite{K01} orthogonality alone of a set $\Lambda$ can be decided by looking at the difference set and making sure that it is contained in the zero set of the short time Fourier transform given by: 
$$
\cZ(V_g g) = \{(t,\nu): V_gg(t,\nu)=0\}.
$$

\begin{cor} \label{orth}
Let $g\in L^2(\R^d)$ such that $\|g\|_2=1$.  
The following statements are equivalent:
\begin{enumerate}
\item[(i)] The Gabor system $\{g_{\la}\}_{\la \in \La}$ is orthonormal;
\item[(ii)] $\La-\La \subseteq \cZ (V_g g ) \cup \{0\}$.
\end{enumerate}
\end{cor}

\subsection{Packing regions}

By Theorem \ref{Gonb char}, for $\|g\|_2=1$, the Gabor system $\{g_\la\}_{\la \in \La}$ is orthonormal if and only if
$$
|V_gg|^2 +\La\leq 1
$$
In order to determine when this system is also complete, and so $\{g_\la\}$ is an orthonormal basis, 
we have to decide whether or not the equality
$$
|V_gg|^2 +\La=1
$$
holds. In general, it is easier to prove that a characteristic function tiles with some set $\La$ than to prove that $|V_gg|^2$ tiles with $\La$. The following simple lemma enables us to check the tiling properties of $|V_gg|^2$ by proving the tiling properties of some characteristic functions:

\begin{lem}\label{glw}
Let $F, G\in  L^1(\R^d)$ be two functions so that  $F, G \geq 0$ and 
$$
\int_{\R^d} F(x)dx=\int_{\R^d} G(x)dx=1.
$$  
Let $\mu $ be a positive Borel measure on $\R^d$ so that $F * \mu\leq 1$  and $G *\mu \leq1$. 
Then, $F * \mu = 1$ if and only if $G * \mu = 1$.
\end{lem}

So, we can replace $|V_gg|^2$ by a simpler function. The proof of this result, as well as the next definition, can be found for instance in  \cite{GLW} (see also the references therein). 

\begin{fed}\label{packreg}
Given $g\in L^2(\R^d)$,  a region $D\subset \R^{2d}$ is an \textbf{(orthogonal) packing region} for $g$  if 
$$
(D^{\circ}-D^{\circ}) \cap \cZ(V_g g)=\varnothing, 
$$
where $D^{\circ}$ denotes the interior of $D$. If $\|g\|_2=1$ then the orthogonal packing region $D$ is called \textbf{tight} when $|D|=1$. 
\end{fed}

\medskip
A simple computation shows that if $D$ is a packing region for $g$ then 
$$
D+\La\leq \R^d
$$
for any $\La$ such that $\La-\La\subseteq \cZ(V_g g)$. Hence, by Lemma \ref{glw}, in order to prove that $|V_gg|^2$ tiles $\R^{2d}$ with $\La$ it is enough to find a tight packing region and prove that $D$ tiles $\R^{2d}$ with $\La$. More precisely

\begin{teo}\label{tiling and packing}
Let $g$ be a norm-one element of $L^2(\R^d)$, $D$ is a tight packing region for $g$, and $\La$ a discrete set such that $\{g_\la\}_{\la\in\La}$ is an orthonormal system. Then the following statements are equivalent:
\begin{itemize}
\item[i.)] The system $\{g_\la\}_{\la\in\La}$ is an orthonormal basis of $L^2(\R^d)$
\item[ii.)] The tight packing region $D$  tiles $\R^{2d}$ with $\La$.
\end{itemize} 
\end{teo}

\medskip

Although to prove the tiling properties of $\chi_D$ is easier than to prove the tiling properties of $|V_gg|^2$, the issue now is to find such a tight packing region when it exists. 

\section{The case of two intervals}

In this section we study Gabor orthonormal bases generated by the characteristic function $g=\chi_\W$, where $\W$ is the union of two intervals. The orthonormality condition imposes that $|\W|=1$, and using the symmetries of the problem, we can restrict our attention to sets of the form
\begin{equation}\label{eq forma del los intervalos}
\W=[0,\alpha) \cup [\alpha+\beta,1+\beta),\quad \quad\mbox{with $\alpha\in [0,1/2]$, and $\beta>0$}.
\end{equation}
If $\alpha<1/2$ we characterize all the possible values of $\beta$, such that for some $\La$, the family $\{g_{\la} \}_{\la \in \La}$ forms an orthonormal basis. In particular, we prove that there exists $\La$ such that $\{g_{\la} \}_{\la \in \La}$ forms an orthonormal basis if and only if $\W$ tiles $\R$ (see Theorem \ref{caso menor}). Moreover, we show that in this case the set $\La$ has a very precise structure (see Theorem  \ref{necstandard}). If $\alpha=1/2$, we also give a characterization of all the possible values of $\beta$, but this time assuming that $\beta\geq 1/2$ (see also Theorem \ref{caso menor}). 
The case $\alpha=1/2$ and $\beta< 1/2$ remains open. To begin with, we present the results, and we leave the proofs for the next two subsections.


\medskip

Our first goal is to provide a geometric characterization of the Gabor orthonormal basis $\{g_{\la} \}_{\la \in \La}$ when $g=\chi_\W$ and $\W$ has the form shown in \eqref{eq forma del los intervalos}. With this aim, we prove that $\chi_\W$ admits a tight orthogonal packing region provided $\alpha<1/2$ or ($\alpha=1/2$ and $\beta\geq 1/2$).

\begin{thm} \label{geom cond 1}
Let $\W=[0,\alpha) \cup [\alpha+\beta,1+\beta)$, where $\alpha\in (0,1/2)$ and $\beta>0$.
Then, the set
$$
D=\W\times [0,1)=\Big([0,\alpha) \cup [\alpha+\beta,1+\beta)\Big)\times [0,1)
$$
 forms a tight orthogonal packing region for $\chi_\W$.
\end{thm}

When both intervals of $\W$ have the same length, the structure of the zero set of $V_{\chi_\W} \chi_\W$ is more complicated. So the packing region is more complicated too.

\begin{thm} \label{geom cond 2}
Let $\W=[0,1/2) \cup [1/2+\beta, \beta+1)$ so that  $\beta\geq 1/2$.
Then, a tight orthogonal packing region for $\chi_\W$ is 
$$
D=\begin{cases}
\displaystyle\W\times  \left( \bigcup_{k=0}^{2\beta} \left[\frac{2k}{2\beta+1}, \frac{2k+1}{2\beta+1}\right )\right)&\mbox{if $\beta\in\frac{1}{2}\N$}\\
\displaystyle\W\times  \left( \bigcup_{k=0}^{ \lfloor 2 \beta \rfloor} \left[\frac{2k}{2\beta+1}, \frac{2k+1}{2\beta+1} \right ) 
\cup \left[ \frac{2(\lfloor 2\beta\rfloor +1)}{2\beta+1}, \frac{2(\lfloor 2\beta\rfloor +1) +\{2 \beta \}}{2\beta+1}\right) \right)&\mbox{if $\beta\notin\frac{1}{2}\N$}
\end{cases}
$$
where $\lfloor x\rfloor$ denotes the (floor) integer part of $x$ and $\{x\}=x-\lfloor x\rfloor$.
\end{thm}

Note that in order to construct $D$, we first consider the set $\W \times \big[0,\frac{1}{2\beta+1}\big)$. 
 Then, we translate it along the axis of $y$ by $\frac{2k}{2\beta+1}$ where $k=1,\dots,2\beta$ if $\beta\in\frac{1}{2}\N$. The set $D$ is just the union of $\W \times \big[0,\frac{1}{2\beta+1}\big)$ and its translates  (see next picture).
 
 \begin{figure}[H]
           \centering
           \includegraphics[height=7cm]{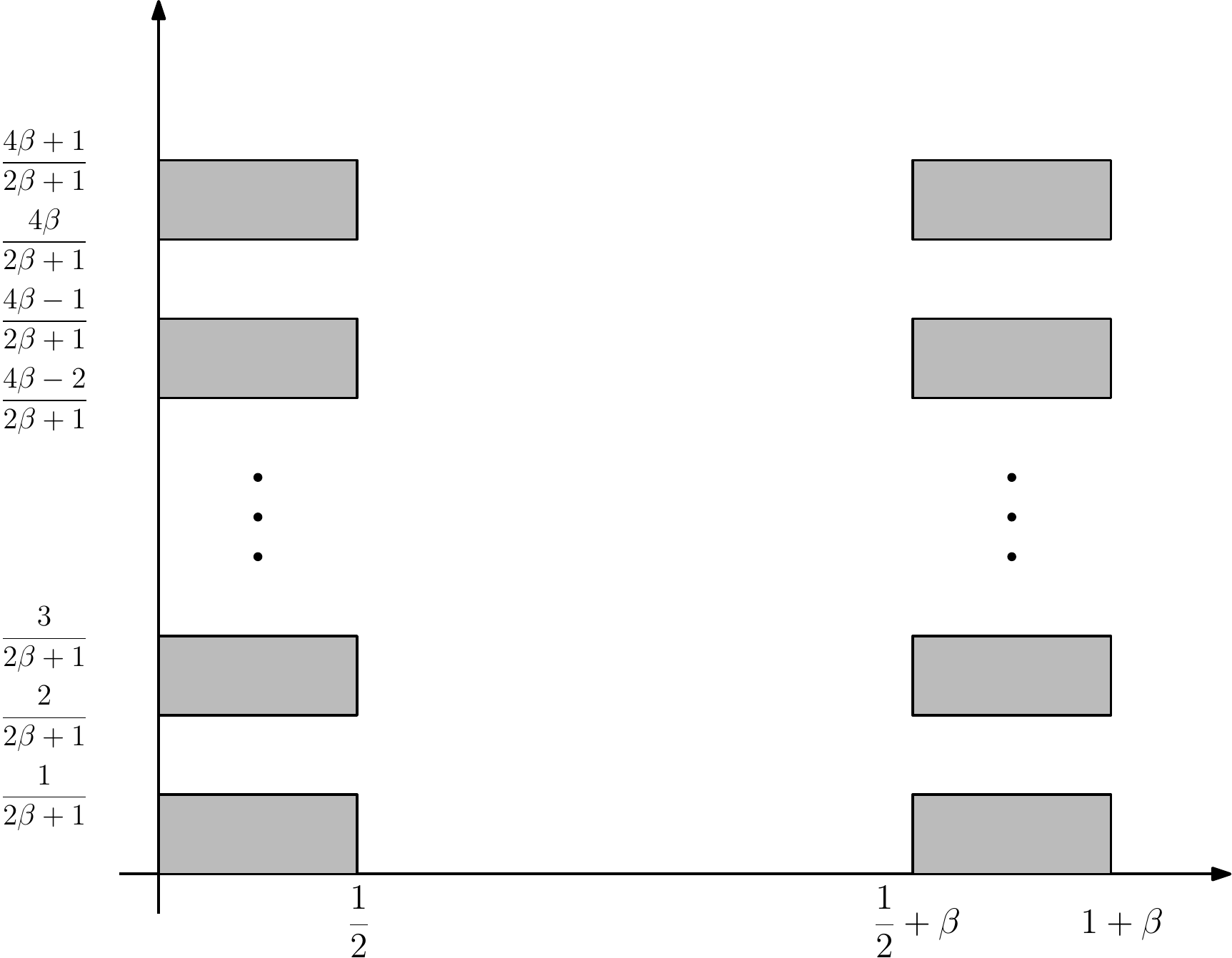}
           \par\vspace{0cm}
           \caption{The set $D$ when  $\beta \in \frac{\N}{2}$}
           \label{Dhalf}
\end{figure}

In the case $\beta\notin\frac{1}{2}\N$, we proceed in a similar way: we consider the set $\W \times \big[0,\frac{1}{2\beta+1}\big)$ and translate it by 
$\frac{2k}{2\beta+1}$ where  $k=1,\dots,\lfloor 2 \beta \rfloor$. However, the union of these sets is not of measure 1. So,  in order to form $D$ so that $|D|=1$, we have to consider the union of the previous sets and we should also add the set 
 $$
 \W\times \left[ \frac{2(\lfloor 2\beta\rfloor +1)}{2\beta+1}, \frac{2(\lfloor 2\beta\rfloor +1) +\{2 \beta \}}{2\beta+1}\right).
 $$
This last set is the union of the thinner sets in the next picture.

 \begin{figure}[H]
           \centering
           \includegraphics[height=7cm]{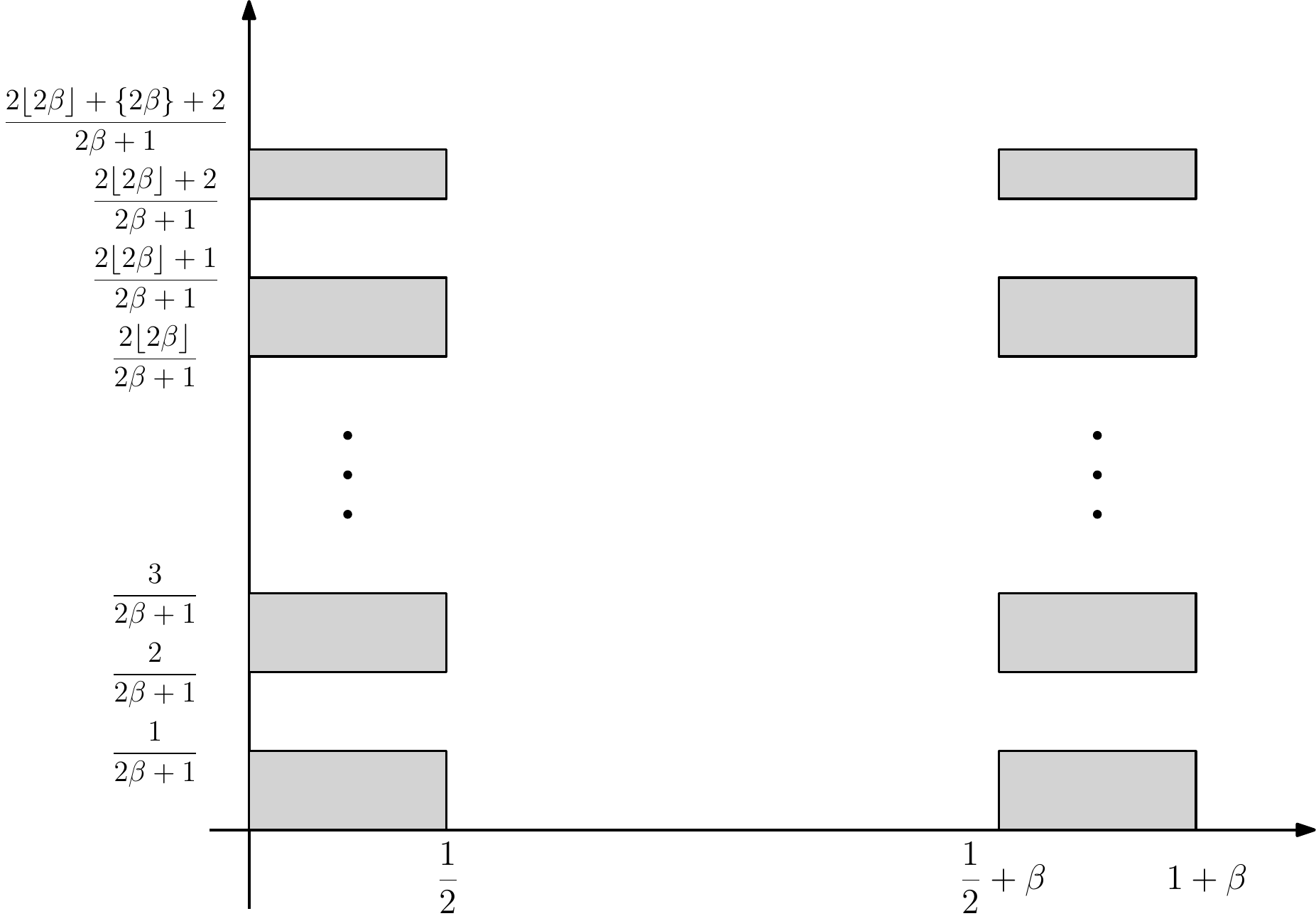}
           \par\vspace{0cm}
           \caption{The set $D$ when $\beta \notin \frac{\N}{2}$}
           \label{Dhalf}
\end{figure}

\medskip

The proofs of  Theorems \ref{geom cond 1} and \ref{geom cond 2} are left for the subsections \ref{section case menor} and \ref{section case igual}. Now, by Corollary \ref{orth} and Theorem \ref{tiling and packing}, Theorems \ref{geom cond 1} and \ref{geom cond 2} directly lead  to the following geometric characterization of the sets $\La$ that produce Gabor orthonormal basis, which is the main result of this paper.

\begin{teo} \label{geom cond}
Let $g=\chi_{\W}$, where $\W=[0,\alpha) \cup [\beta+\alpha, \beta+1)$, where $\alpha<1/2$ or ($\alpha=1/2$ and $\beta\geq 1/2$).
Then $\{g_{\la}\}_{\la \in \La}$ forms an orthonormal basis if and only if
\begin{itemize}
\item[(i)] $\La-\La \subset \cZ (V_g g) \cup \{0\}$, 
\item[(ii)] $\La+ D=\R^2$,
\end{itemize}
where $D$ is the tight othogonal packing region described by Theorem \ref{geom cond 1} or Theorem \ref{geom cond 2}  depending on the case.
\end{teo}

\medskip

\begin{rem}
When $\alpha=1/2$ and $\beta<1/2$, the structure of the zero set of $V_{\chi_\W} \chi_\W$ is even more complicated. From the computations that we will show in Subsection \ref{section case igual} it follows that 
$$
D=[0,1+\beta)\times[0,(1+2\beta)^{-1})
$$
is a packing region for $\chi_\W$. However $|D|<1$. We could not prove the existence of a \textbf{tight} packing region. So, our techniques can not be applied in this case.  \EOE
\end{rem}

\subsection*{Gabor orthonormal basis and tiling properties of $\W$}

Our next goal is to relate the existence of Gabor orthonormal basis having as window the function $\chi_\W$ with some tiling properties of $\W$. In \cite{L}, $\L$aba proved that the Fuglede conjecture is true for the union of two intervals. Thus, if a set $\W$ tiles the real line with some set $\Gamma$, then it is also spectral, which means that $L^2(\W)$ admits an orthonormal basis of exponentials $\{e_\la\}_{\la\in\La}$.  More precisely she got the next result:

\begin{teo}\label{LabaT}
Let $\W$ be a union of two intervals as in equation \eqref{eq forma del los intervalos}. Then $\W$ tiles the real line if and only if  it is spectral. Moreover, this happens  either if $\alpha<1/2$ and $\beta\in\Z$ or if $\alpha=1/2$ and $\beta\in \frac12 \Z$.
\end{teo}

\medskip

As a consequence of Theorem \ref{geom cond} and Fubini's theorem we get the following result.


\begin{teo}\label{caso menor}
Let $\W=[0,\alpha) \cup [\alpha+\beta,1+\beta)$, where $\alpha<1/2$ or ($\alpha=1/2$ and $\beta\geq 1/2$). Then, there exists a Gabor orthonormal basis associated to the function $\chi_\W$ if and only if $\W$ tiles the real line. This occurs exactly when any of the following conditions holds:
\begin{itemize}
\item[(i)]  $\displaystyle 0<\alpha < 1/2 $ and  $\beta \in \N$; 
\item[(ii)] $\alpha =1/2$  and $\displaystyle \beta \in \frac{1}{2}\N$.
\end{itemize} 
\end{teo}

\bdem  
Let $g=\chi_{\W}$. Let $\La$ be a subset of $\R^2$ so that $\{g_{\la}\}_{\la\in \La}$ is a Gabor orthonormal basis. 
By Theorem \ref{geom cond} we get that $ D+\La=\R^2$, where $D$ is a tight orthogonal packing region 
(explicitly given in Theorems \ref{geom cond 1} and \ref{geom cond 2}). Since $D$, which is a Cartesian product of $\W$ with another set, tiles $\R^2$, by Fubini's theorem $\W$ must tile the line. The description of such sets is given in \cite{L}.
The converse is easier, because if $\W$ tiles the real line,  it also admits a spectrum (see  \cite{L}).  Therefore, as we mentioned in the introduction, there exists a set $\La$ so that $\{g_{\la}\}_{\la\in \La}$ is a Gabor orthonormal basis.
\edem

\medskip

Let $\displaystyle  g = \chi_{\W}$, where $\W\subset \R^d$ is measurable of measure 1. Suppose that 
$$
\La= \bigcup_{t\in J} \{t\} \times \La_t
$$
where
\begin{itemize}
\item[(i)] $\bigcup_{t\in J} (\W+t)= \R^d$, and
\item[(ii)] For every $t\in J$ the system $\{e^{2\pi i \pint{\la,x}}\}_{\la \in \La_t}$ is an orthonormal basis for $L^2(\W+t)$, which is equivalent to being an orthonormal basis for  $L^2(\W)$.
\end{itemize} 

When $\La$ has this structure it is called \textbf{$\W$-standard}. It is not difficult to prove that in this case $\{g_{\la}\}_{\la\in \La}$ is an orthonormal basis for $L^2(\R^d)$. The opposite implication, that is, the necessity of $\La$ to be $\W$-standard when $\{g_{\la}\}_{\la\in \La}$ is a basis, is not always true. Moreover, this is not necessary even in the case when $g= \chi_{[0,1)^d}$ and  $d\geq 2$. In fact, in \cite{GLW} it was proved that in this case there exist sets $\La$  so that $\{g_{\la}\}_{\la \in \La}$ forms a Gabor orthonormal basis, but the sets $[0,1)^d+\la$ with  $\la\in\La$ have significant overlaps.

\medskip

 However, in the same paper it was shown that, if $d=1$ and so $g=\chi_{[0,1)}$, then the system $\{g_{\la}\}_{\la \in \La}$ is a Gabor orthonormal basis of $L^2(\R)$ if and only if $\La$ is standard. The same holds in our case, that is, when  $\W$ has the form described by Theorem \ref{caso menor}.

\medskip

\begin{teo} \label{necstandard}
Let $g=\chi_{\W}$, where $\W=[0,\alpha) \cup [\alpha+\beta,1+\beta)$ and $\alpha<1/2$. If the system $\{g_{\la}\}_{\la \in \La}$ is a Gabor orthonormal basis for $L^2(\R)$, then $\La$ is standard, i.e.
$$
\La=\bigcup_{k\in\Z}\{k\}\times (\Z+a_k),
$$
where $a_k\in[0,1)$ for every $k\in\Z$.
\end{teo}

\begin{teo} \label{necstandard2}
Let $g=\chi_{\W}$, where $\W=[0,\alpha) \cup [\alpha+\beta,1+\beta)$, $\alpha=1/2$ and $\beta\in \frac{\N}{2}$. If the system $\{g_{\la}\}_{\la \in \La}$ is a Gabor orthonormal basis for $L^2(\R)$, then $\La$ is standard, i.e.
\begin{equation}\label{expl}
\La=\bigcup_{k\in K}\left\{\frac{k}{2}\right\}\times\left(\frac{L_k+a_k}{2\beta+1}\right),
\end{equation}
where
\begin{itemize}
\item[i.)] $K, L_k\subseteq \Z$ for every $k\in K$;
\item[ii.)] $a_k\in [0,1)$ for every $k\in K$;
\item[iii.)] $K\cup ((2\beta+1)+K)=\Z$;
\item[iv.)] $L_k+\{2n:\ n=0,\pm1,\ldots, 2\beta\}=\Z$.
\end{itemize}
\end{teo}

\medskip

The proofs of these result are given in the next two subsections. We will consider separately the case $\alpha<1/2$ and the case $\alpha=1/2$.

\subsection{The case $\alpha<1/2$}\label{section case menor} 

This section contains the proofs of Theorems \ref{geom cond 1} and \ref{necstandard}. 
We will start this section with some technical lemmata.

\begin{lem}\label{symmetries}
Let $\W\subseteq \R^d$. If $g=\chi_\W$, then 
$$
|V_gg(t,\nu)|=|V_gg(t,-\nu)|=|V_gg(-t,\nu)|.
$$
\end{lem}
\bdem
Since $g$ is a real valued function, $V_gg(t,-\nu)=\overline{V_gg(t,\nu)}$.  
This proves the first equality. For the second one, observe that
$$
|V_gg(t,\nu)|=|\widehat{\chi}_{\W\cap(\W+t)}(\nu)|=|e^{-2\pi i \nu t}\widehat{\chi}_{\W\cap(\W-t)}(\nu)|=|V_gg(-t,\nu)|.
$$
\edem

\begin{lem}\label{easy}
If $I$ is a bounded interval, then $\widehat{\chi}_I(\w)\neq 0$ for every $\w\in(-|I|^{-1},|I|^{-1})$.
\end{lem}

\bdem
Indeed, if $2r=|I|$ then it holds that
$$
|\widehat{\chi}_I(\w)|=|\widehat{\chi}_{[-r,r]}(\w)|=\left|\frac{\sin(2\pi r \w)}{\pi \w}\right|.
$$
\edem

\begin{lem}\label{dif} 
Let   $I$ and $J$ be two disjoint intervals, satisfying $|I|<|J|$ and  $|I|+|J| <1$. Then, if 
 $$
 \Omega= I \cup J,
 $$
we have that $\widehat{\chi}_\W(\w)\neq 0$ for every $\w\in[-1,1]$.
\end{lem}

\bdem
We can consider $\w\neq 0$, since clearly $\widehat{\chi_{\W}}(0)=|\W|\neq 0$.
Note that given $a,b\in \R$,  we have 
\begin{equation}\label{patata frita}
\widehat{\chi_{[a,b]}}(\w) = \frac{\sin \pi \w (b-a)}{\pi \w}  \, \, e^{-\pi i (a+b) \,\w}.
\end{equation}
Let $\ell_1=|I|$, $\ell_2=|J|$, and  let $m_1$,  $m_2$ be the midpoints of $I$ and $J$ respectively. Then, if 
\begin{align*}
\widehat{\chi_{\W}}(\w)
&= e^{-2\pi i m_1 \w}   \frac{\sin \pi \ell_1 \w}{\pi \w}
+ e^{-2 \pi i  m_2  \w } \frac{\sin \pi \ell_2 \w}{\pi \w}=0 
\end{align*}
we get that
$$
| \sin (\pi \ell_1  \w) | =|\sin (\pi \ell_2 \w)|.
$$
This is not possible for $|\w|\leq 1$,  since $0<\ell_1 <\ell_2$ and $\ell_1+\ell_2<1$. Hence, 
$\widehat{\chi_{\W}}(\w)\neq 0$ for $\w\in[-1,1]$.
\edem

Finally, when $\beta<\alpha$, then the intersection of the original set with the same set translated by $t\in[\beta,\alpha)$ gives a union of three intervals,  illustrated by the following picture:

 \begin{figure}[H]
           \centering
           \includegraphics[height=4cm]{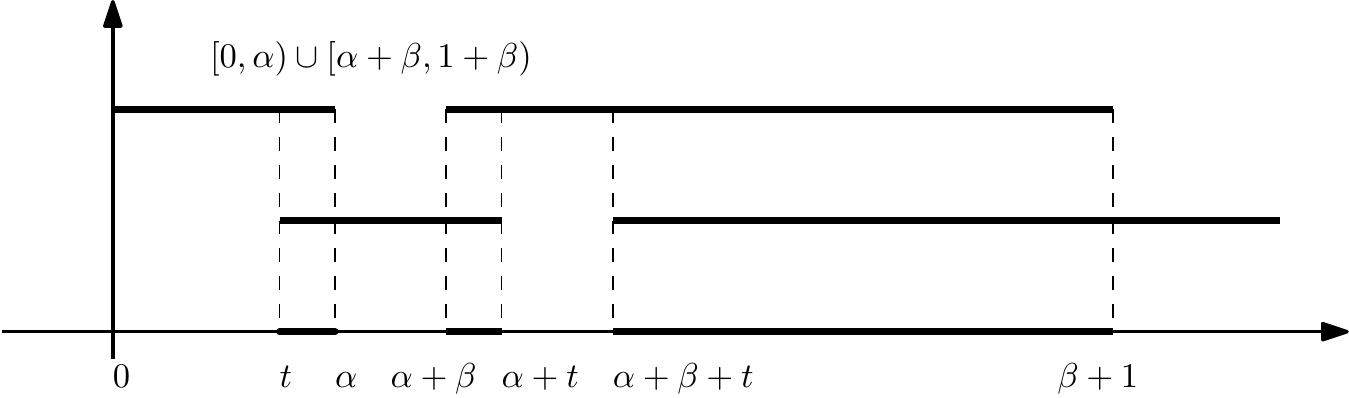}
           \par\vspace{0cm}
           \caption{Example}
           \label{Dhalf}
\end{figure}

\begin{lem}\label{dif3}
Given $0<\beta<t<\alpha<1/2$, if 
$$
\W =[t,\alpha)\cup [\alpha+\beta,\alpha+t)\cup [\alpha+\beta+t,1+\beta)
$$
then $\widehat{\chi}_\W(\w)\neq 0$ for every $\w\in[-1,1]$. 
\end{lem}

\bdem First of all, note that
$$
\chi_\W=\chi_{[t,1+\beta)}-\big(\,\chi_{[\alpha\,,\,\alpha+\beta)}+\chi_{[\alpha+t\,,\,\alpha+\beta+t)}\,\big).
$$
Therefore, using \eqref{patata frita} we get for $\w\neq 0$
\begin{align*}
\widehat{\chi}_\W(\w)&=e^{-i\pi  (1+t+\beta)\w}\ \frac{\sin(\pi \w (1+\beta+t))}{\pi \w }\,-\,
 \frac{\sin(\pi \w \beta)}{\pi \w}\Big(e^{-i\pi \w(2\alpha+\beta)}+e^{-i\pi \w(2\alpha+2t+\beta)}\Big) \\
&=e^{-i\pi  (1+t+\beta)\w} \left(\frac{\sin(\pi \w (1+\beta+t))}{\pi \w }- \frac{\sin(\pi \w \beta)}{\pi \w}\Big(e^{i\pi \w(1-2\alpha+t)}+e^{i\pi \w(1-2\alpha-t)}\Big)\right).
\end{align*}
Suppose there exists  $|\w|\leq 1$ such that $\widehat{\chi}_\W(\w)=0$. Clearly $\w\neq 0$. So, we get
$$
\frac{\sin(\pi \w (1+\beta+t))}{\pi \w }= \frac{\sin(\pi \w \beta)}{\pi \w}\Big(e^{i\pi \w(1-2\alpha+t)}+e^{i\pi \w(1-2\alpha-t)}\Big).
$$
Since $\beta\in(0,1/2)$, comparing the imaginary parts of both sides we obtain that 
\begin{equation}\label{arkuda}
\sin(\pi\w(1-2\alpha+t))=-\sin(\pi \w(1-2\alpha-t)).
\end{equation}
Recall that $0<t<\alpha<1/2$, therefore
$$
-\frac{1}{2}<|\w| (1-2\alpha-t)<|\w|(1-2\alpha+t)<1-\alpha<1.
$$
In consequence, the identity \eqref{arkuda} holds if and only if any of the following holds
$$
\pi|\w|(1-2\alpha+t)=-(\pi |\w|(1-2\alpha-t))
$$
or
$$
\pi - \pi|\w|(1-2\alpha+t)=-\pi |\w|(1-2\alpha-t)).
$$
In the first case, $1-2\alpha=0$ which is impossible. In the second case, $2|\w|t=1$ which is also impossible because $|2t|<1$ and $|\w|\leq 1$. This completes the proof. 
\edem

\bdem[Proof of Theorem \ref{geom cond 1}]

We split the proof in two cases:

\medskip

\textbf{(i) Case $\beta \geq \alpha$:} Since $|D|=1$, it is enough to prove that
$(\mathring{D}-\mathring{D})\cap  \cZ (V_g g)=\varnothing$. By the symmetries of $V_gg$ proved in Lemma \ref{symmetries}, 
we only have to prove that the set
$$
\big\{(|x_1-y_1|,|x_2-y_2|):\ x,y\in \mathring{D}\big\}=\Big([0,1-\alpha) \cup (\beta,\beta+1)\Big)\times [0,1)
$$
does not intersect $\cZ (V_g g)$. With this aim, first note that for $t> 0$
\begin{equation} \label{intersections}
\W\cap(\W+t)=\begin{cases}
[t,\alpha]\cup[\beta+\alpha+t,\beta +1]&\mbox{if $t\in[0,\alpha)$}\\
[\beta+\alpha+t,\beta+1]&\mbox{if $t\in[\alpha,1-\alpha)$}\\
[\beta+\alpha,t+\alpha]&\mbox{if $t\in[\beta,\beta+\alpha)$}\\
[t,t+\alpha]&\mbox{if $t\in[\beta+\alpha,\beta+1-\alpha)$}\\
[t,\beta+1]&\mbox{if $t\in[\beta+1-\alpha,\beta+1)$}
\end{cases}
\end{equation}
and is empty otherwise.
Therefore, using Lemmas \ref{easy} and \ref{dif} we get that $\widehat{\chi}_{\W\cap(\W+t)} (\nu)\neq 0$, for every $\nu \in [0,1)$.  
Since $ V_gg(t,\nu)=\widehat{\chi}_{\W\cap(\W+t)}(\nu) $, we have that
$$
\big([0,1-\alpha) \cup (\beta,\beta+1)\big)\times [0,1)
$$
is  region free of zeroes of $V_gg$. \\
 
\textbf{ Case $0<\beta<\alpha$}: The idea of the proof is exactly the same, but in the description of $\W\cap(\W+t)$ appears a new case, when $t\in[\beta,\alpha)$:
\begin{equation} \label{intersectionshalf}
\W\cap(\W+t)=\begin{cases}
[t,\alpha] \cup[\beta+\alpha+t,\beta +1]&\mbox{if $t\in[0,\beta)$}\\
[t,\alpha) \cup [\beta+\alpha,\alpha +t) \cup [\alpha+\beta+t, \beta+1)&\mbox{if $t\in[\beta,\alpha)$}\\
[\beta+\alpha,\alpha+t) \cup [\alpha+\beta+t,\beta+1)&\mbox{if $t\in[\alpha,\beta+\alpha)$}\\
[t,\alpha+t)\cup [\alpha+\beta+t, \beta+1 )&\mbox{if $t\in[\beta+\alpha,\beta+1-\alpha)$}\\
[t,\beta+1]&\mbox{if $t\in[\beta+1-\alpha, \beta+1)$}
\end{cases}
\end{equation}
and is empty otherwise. If $ t\in[\beta,\alpha)$ we get that $\widehat{\chi}_{\W\cap(\W+t)} (\nu)\neq 0$ by Lemma \ref{dif3}.
For the rest of the cases,  as previously we use Lemmas \ref{easy} and \ref{dif} to get that $\widehat{\chi}_{\W\cap(\W+t)} (\nu)\neq 0$, 
for every $\nu \in [0,1)$.  Since $ V_gg(t,\nu)=\widehat{\chi}_{\W\cap(\W+t)}(\nu) $, we have that
$$
\big([0,1-\alpha) \cup (\beta,\beta+1)\big)\times [0,1)
$$
is  region free of zeroes of $V_gg$. 
\edem

Now we will proceed to the proof of Theorem \ref{necstandard}. We will start with the following definition.

\begin{defi}\label{Tannhauser}
A pair $(A,B)$ of bounded open sets in $\R^d$ is called a  \textbf{ spectral pair} if
$$
(A^{\circ}-A^{\circ}) \cap \cZ(\widehat{\chi}_{B})=\varnothing  
\quad \text{ and } \quad
(B^{\circ}-B^{\circ}) \cap \cZ(\widehat{\chi}_{A})=\varnothing.
$$
If $|A|=|B|=1$ we say that the spectral pair is \textbf{tight}.
\end{defi}

\begin{exas}\label{examples de SP} 
$\ $
\begin{itemize} 
\item[(i)] The pair  $([0,1), [0,1))$ is a tight spectral pair. 

\item[(ii)] Let  $\W= [0,\alpha) \cup [\alpha+\beta,1+\beta)$, where $\displaystyle 0<\alpha < 1/2 $ and $\beta \in \N$. Then, the pair  $(\W, [0,1))$ is tight spectral. In fact, by Lemma \ref{easy} we get
$
 \widehat{\chi_{[0,1)}}(\w)  \neq 0
$ 
for every $\w \in (-1,1)$ and  since  
$$
(\W^{\circ}-\W^{\circ}) =(\alpha-1,1-\alpha)\cup(\beta,\beta+1)\cup(-(\beta+1),-\beta),
$$ 
we have that 
$$
(\W^{\circ}-\W^{\circ}) \cap \cZ(\widehat{\chi_{[0,1)}})=\varnothing.
$$
On the other hand, by Lemma \ref{dif} we know that
$
\widehat{\chi_{\W}}(\w) \neq 0
$ 
for every $\w \in [-1,1]$. Therefore
$$
([0,1)^{\circ}-[0,1)^{\circ})  \cap \cZ(\widehat{\chi_{\W}})=\varnothing.
$$
\end{itemize}\EOE
\end{exas}

\begin{lem}\label{prod}
Let $(A, B)$ and $(C,D)$ be (tight) spectral pairs of bounded open sets $A, B,C,D$ in $\R^d$. Then, the pair 
$(A\times C, B\times D)$ is a (tight) spectral pair. 
\end{lem}

The proof is easy and is ommitted.

\begin{rem}\label{ejemplos de productos}
By Lemma \ref{prod} the pair $(\W\times [0,1), [0,1)^2)$ is also a  tight spectral pair. \EOE
\end{rem}

\medskip

\noi Next result was proved in \cite{K00} (see Theorem 9).

\begin{teo}\label{KOO}
Assume that $(A,B)$ is a tight spectral pair. Then $\La$ is a spectrum of $A$ if and only if
$B+ \La$ is a tiling.
\end{teo}

Note that as a consequence of this theorem, and the above mentioned examples we get the following corollary.

\begin{cor}\label{Spec vs tilign en cubos}
A set $\La$ is a spectrum of $[0,1)^d$ if and only if $[0,1)^d+ \La$ is a tiling of $\R^d$.
\end{cor}


Now, we are ready to prove why $\La$ should be standard. The proof follows similar lines as in \cite{GLW}, 
but for the sake of completeness we present it here.

\bdem [Proof of Theorem \ref{necstandard}]
Let $g=\chi_{\W}$. If the system $\{g_{\la}\}_{\la \in \La}$ 
is a Gabor orthonormal basis of $L^2(\R)$,  by Theorem \ref{geom cond} we get that $D+\La= \R^2$, where $D=\W\times [0,1)$.
If $\Delta=[0,1)^2$, then by Remark \ref{ejemplos de productos}  $(D, \Delta)$ forms a tight spectral pair. So,  by Theorem \ref{KOO} 
$\La$ is a spectrum for $\Delta$.  By Corollary \ref{Spec vs tilign en cubos}, $\La+\Delta$ is a tiling of $\R^2$. Therefore, by simple inspection (a tiling of the plane by a square is either by \lq \lq shifted columns" or by \lq\lq shifted rows"),  the set $\La$ can be of any of the following form:
\begin{equation} \label{ins}
\La=\bigcup_{k\in \Z} (\Z+a_k)\times \{k\} \quad \text{ or } \quad \La=\bigcup_{k\in \Z} \{k\} \times  (\Z+a_k),  \
\end{equation}
where $a_k$ are real numbers in $[0,1)$ for $k\neq 0$ and $a_0=0$. It only remains to prove that $\La$ cannot be of the form 
\begin{equation} \label{insl}
\La=\bigcup_{k\in \Z} (\Z+a_k)\times \{k\}, 
\end{equation}
unless $a_k=0$ for every $k\in\Z$. Assume that there exists $a_k\neq 0$. By symmetry, we can assume that $k=\min\{n\in\N: a_n\neq 0 \}>0$. Since  $\beta \in \Z$ and by definition of $k$
$$
\la=(a_k+\beta, k)\in \La \quad \text{ and } \quad  \mu=(0, k-1)  \in \La.
$$
Hence $\la-\mu=(a_k+\beta, 1)$ should belong in $\cZ(V_g g)$ by (i) of Theorem \ref{geom cond}. 
However we will see that the point $(a_k+\beta, 1)$ belongs to a region free of zeroes of $V_gg$. 
Thus, we get a contradiction and so $\La$ cannot be as in \eqref{insl}.
It remains to see that the point $(a_k+\beta, 1)$ belongs to a region free of zeroes of $V_gg$. Note that
 $t:=\alpha_k +\beta \in (\beta,\beta+1)$. By Theorem \ref{caso menor}, $\beta\in \N$ and so  $\beta >\alpha$.
Since 
$$
V_g g(t, \nu) =\widehat{\chi}_{\W \cap (\W+t)} (\nu)
$$
for every $t\in (\beta,\beta+1)$ by \eqref{intersections} we get
\begin{equation*} 
\W\cap(\W+t)=\begin{cases}
[\beta+\alpha,t+\alpha]&\mbox{if $t\in[\beta,\beta+\alpha)$}\\
[t,t+\alpha]&\mbox{if $t\in[\beta+\alpha,\beta+1-\alpha)$}\\
[t,\beta+1]&\mbox{if $t\in[\beta+1-\alpha,\beta+1)$.}
\end{cases}
\end{equation*}
Therefore, using Lemma \ref{easy}  we get that $ V_gg(t,\nu)\neq 0$ and so
$$
 (\beta,\beta+1)\times [0,1]
$$
is  region free of zeroes of $V_gg$. Therefore the point $(a_k+\beta, 1)\not \in \cZ(V_g g)$.
\edem

\subsection{The case $\alpha=1/2$ and  $\beta\geq 1/2$}\label{section case igual}

Now we will prove Theorem \ref {geom cond 2}. The idea is the same, but due to the extra symmetries, the zero set of $V_{\chi_\W}\chi_\W$ are different. 

\begin{lem}\label{dif2}
Let $I$ and $J$ be two disjoint intervals so that $|I|=|J|\leq1/2$ and let $m_1, m_2$, with $m_1 < m_2$, be the midpoints of $I, J$ respectively. Then, if 
$$
\W =I\cup J,
$$
and $\w\in(-2,2)$, then  $\widehat{\chi}_\W(\w)=0$ if and only if $\w=\frac{k}{2(m_2-m_1)}$ 
for a non zero odd integer $k$ such that $|k|<4(m_2-m_1)$. 
\end{lem}

\bdem 
Since $\widehat{\chi}_\W (0)=|\W|\neq 0$, we can consider $\w\neq 0$. Let $\ell=|I|=|J|$. 
Then
\begin{align*}
\widehat{\chi_{\W}}(\w)
&=\frac{\sin \pi \ell \w}{\pi \w}   ( e^{-2\pi i m_1\w} + e^{-2\pi i m_2\w }).
\end{align*}
On the one hand, since $\ell\leq1/2$ and $\w \in (-2,2)\setminus \{0\}$, the term 
$$
\frac{\sin \pi \ell \w}{\pi \w} \neq 0.
$$
On the other hand 
$$
 e^{-2\pi i m_1\w} + e^{-2\pi i m_2\w }=0
$$
if and only if  $\w= k/2(m_2-m_1)$ for an odd integer $k$. Moreover, since $\w\in (-2,2)\setminus \{0\}$ we have that $|k|<4(m_2-m_1)$.
\edem

\bdem[Proof of Theorem \ref{geom cond 2}]


\textbf{(i) Case $\displaystyle \beta \in \frac{1}{2}\N $:} Note that in this case, for $t>0$
\begin{equation} \label{intersectionshalf}
\W\cap(\W+t)=\begin{cases}
[t,1/2] \cup[\beta+1/2+t,\beta +1]&\mbox{if $t\in[0,1/2)$}\\
[\beta+1/2,t+1/2]&\mbox{if $t\in[\beta,\beta+1/2)$}\\
[t,\beta+1]&\mbox{if $t\in[\beta+1/2,\beta+1)$}
\end{cases}
\end{equation}
and is empty otherwise.
Let
\begin{align*}
D&= \W\times  \left( \bigcup_{k=0}^{2\beta} \left[\frac{2k}{2\beta+1}, \frac{2k+1}{2\beta+1}\right )\right).
\end{align*}
Note that 
$$
|D|= \sum_{0}^{2\beta}\frac{1}{2\beta+1}=1.
$$
Hence, as in the case of Theorem \ref{geom cond 1}, by the symmetry of  $ V_gg(t,\nu)$ it is enough to prove that 
$$
\big\{(|x_1-y_1|,|x_2-y_2|):\ x,y\in \mathring{D}\big\}
$$
does not intersect the zero set of $ V_gg(t,\nu)=\widehat{\chi}_{\W\cap(\W+t)}(\nu)$.  Given $x,y\in \mathring{D}$, it holds that $|x_2-y_2|\leq 2$. Therefore, by Lemmas \ref{easy} and \ref{dif2}, the unique possibility of intersection is when $|x_1-y_1|$ is in $[0,1/2]$. In this region, by Lemma \ref{dif2} the zeroes are located at the heights 
$$
\frac{k}{2\beta+1}
$$
for $k$ an odd integer such that $k\leq 2(2\beta+1)$. So, we place the blocks in such 
way that we avoid these lines of zeroes  (see figure \ref{Dhalf} above). \\

\textbf{(ii) Case $\displaystyle \beta\not \in \frac{1}{2} \N $.} 
We can construct  $D$ as the following union  
\begin{align*}
D&= \W\times  \left( \bigcup_{k=0}^{ \lfloor 2 \beta \rfloor} \left[\frac{2k}{2\beta+1}, \frac{2k+1}{2\beta+1} \right ) 
\cup \left[ \frac{2(\lfloor 2\beta\rfloor +1)}{2\beta+1}, \frac{2(\lfloor 2\beta\rfloor +1) +\{2 \beta \}}{2\beta+1}\right) \right)\\
\end{align*} 
Note that $ \beta\not \in \frac{1}{2} \N $ and so we have considered the union  until $\lfloor 2 \beta \rfloor$. 
Since $k\leq 2(2\beta+1)$ we have also added the set
$$ 
A=\left[ \frac{2(\lfloor 2\beta\rfloor +1)}{2\beta+1}, \frac{2(\lfloor 2\beta\rfloor +1) +\{2 \beta \}}{2\beta+1}\right) \times \W
$$ 
which is free of zeros. The set $A$ have been added in order to get  $|D|=1$.
\edem

\subsection*{Proof of Theorem \ref{necstandard2}}
 The  proof of this result has been inspired by the proofs of Proposition 3.2 and Theorem 3.3 in \cite{GLW}.

\medskip

To start with, recall that by Theorem \ref{geom cond 2}, a tight orthogonal packing region for $\W$ is:
\begin{align*}
D&=\W\times  \left( \bigcup_{k=0}^{2\beta} \left[\frac{2k}{2\beta+1}, \frac{2k+1}{2\beta+1}\right )\right)\\
&=\Big([0, 1/2) \cup [\beta+1/2, \beta+1)\Big) \times  \left( \bigcup_{k=0}^{2\beta} \Big[\frac{2k}{2\beta+1}, \frac{2k+1}{2\beta+1}\Big )\right).
\end{align*}
If $\{g_{\la}\}_{\la \in \La}$ is a Gabor orthonormal basis of $L^2(\R)$ then, by Theorem \ref{geom cond}, we get that 
$$
D+\La= \R^2.
$$ 
As usual, we will assume that $0\in \La$. Consider the matrix 
\begin{equation} \label{delta}
M:=\begin{pmatrix}
2 & 0\\
0 & 2\beta+1\\
\end{pmatrix},
\end{equation}
and define $D_{M}= M(D)$ and $\La_M=M(\La)$. Clearly $\La_M +D_M = \R^2$, and
\begin{align*}
D_{M}&= \Big([0,1)\cup [2\beta +1, 2\beta+2) \Big) \times   \left( \bigcup_{k=0}^{2\beta} \big[{2k}, {2k+1}\big )\right)\\
&=\bigcup_{\gamma\in\Gamma} \gamma+[0,1)^2,
\end{align*}
where $\Gamma=\{0, 2\beta+1\}\times\{0,2,\ldots,4\beta\}$. Therefore, $(\Gamma+\La_M)+[0,1)^2=\R^2$. By inspection,
this implies that either
\begin{equation} \label{bolero}
\Gamma+\La_M=\bigcup_{k\in \Z} (\Z+a_k)\times \{k\} \quad \text{ or } \quad \Gamma+\La_M=\bigcup_{k\in \Z} \{k\} \times  (\Z+a_k),  \
\end{equation}
where $a_k$ are real numbers in $[0,1)$ for $k\neq 0$, and $a_0=0$ by our initial assumption on $\La$. 

\medskip

Till now, we have used the tiling condition over $D$ in order to get \eqref{bolero}. However, this condition by itself is not enough to prove that 
$\La$ is standard. To achieve this, we have to use the extra structure of the set $\La$, imposed by the orthogonality of the system $\{g_\la\}_{\la\in \La}$. More precisely, we will use that $\La$ satisfies
$$
\La-\La \subset \cZ (V_g g) \cup \{0\}.
$$
To begin with, we prove the following claim.

\paragraph{Claim:} The equality $\displaystyle\Gamma+\La_M=\bigcup_{k\in \Z} (\Z+a_k)\times \{k\}$ does not hold if $a_k\neq0$ for some $k\in\Z$. 

To prove this claim, assume that there exist $k\neq 0$ such that $a_k\neq 0$. By the symmetries of the problem, we can suppose without loss of generality that $k>0$, and that it satisfies the condition
$$
k=\min\{\ell>0: \ a_\ell\neq 0\}.
$$
For each $j\in\Z$, let $\gamma_j,\gamma_j'\in\Gamma$ and $\la_j,\la'_j\in \La_M$ be such that
\begin{equation} \label{contra1}
(a_k+j,k)=\gamma_j+\la_j
\end{equation}
\begin{equation}  \label{contra2}
(j,k-1)=\gamma'_j+\la'_j\ .
\end{equation}
Therefore, for every $j$ we have that
$$
(a_k,1)=(\gamma_j-\gamma'_j)+(\la_j-\la'_j).
$$
Since $\gamma_j-\gamma'_j\in \Gamma-\Gamma$, the possibilities for its first coordinate, denoted by $(\gamma_j-\gamma'_j)_1$, are 
$$
0,\quad  2\beta+1, \peso{and} -(2\beta+1).
$$
Suppose that for some $j$ it holds that  $\gamma_j-\gamma'_j=(2\beta+1,n)$ with 
$$
n\in\{2k:\ k=0,\pm1,\ldots, \pm 2\beta\}.
$$
 Then
$$
(t,\nu):=\la_j-\la'_j= (a_k-(2\beta+1),1-n ).
$$
Since $\{g_{\la}\}_{\la \in \La}$ is a Gabor orthonormal basis of $L^2(\R)$, by Theorem \ref{geom cond} $M^{-1}(t,\nu)\in \cZ (V_g g)$. However, 
$$
|V_gg(M^{-1}(t,\nu))|=|V_gg(t/2,\nu/(2\beta+1))|=|\widehat{\chi}_{\W\cap(t/2+\W)}(\nu/(2\beta+1))|.
$$
Since
$$
|\W\cap(t/2+\W)|=\Big|\,\Big[ \frac{a_k}{2},\frac12\Big)\,\Big|<\frac{1}{2},
$$
and $|\nu/(2\beta+1)|<2$, by Lemma \ref{easy} we get that $|V_gg(M^{-1}(t,\nu))|\neq 0$. This proves that $(\gamma_j-\gamma'_j)_1\neq 2\beta+1$. A similar argument shows that
$(\gamma_j-\gamma'_j)_1\neq -(2\beta+1)$. Moreover, we can also compare the first coordinates of $\gamma_j$ and $\gamma'_{j+1}$, and again the same arguments show that neither $(\gamma_j-\gamma'_{j+1})_1=\pm(2\beta+1)$ is possible. Therefore, we have that 
$$
(\gamma_j-\gamma'_j)_1=(\gamma_j-\gamma'_{j+1})_1=0,\quad \forall j \in \Z.
$$
Equivalently, for the first coordinates, we have
$$
\text{(i)}: (\gamma_j)_1=(\gamma'_{j})_1 \quad \text{and } \quad  \text{(ii)}: (\gamma_{j})_1=(\gamma'_{j+1})_1 \quad \forall j \in \Z.
$$
Fix now $j\in\{0,1,2,\dots,2\beta, 2\beta +1\}$. Then, combining these identities we get that 
$$
 (\gamma_0)_1  \myeqd  (\gamma_1')_1 \myeq (\gamma_1)_1 \myeqd  (\gamma_2')_1  \myeq\dots  \myeqd (\gamma_{2\beta+1})_1 
 \myeq  (\gamma_{2\beta+1}')_1.
$$
Hence,  $\gamma_0-\gamma'_{2\beta+1}=(0,n)$. Note that $n$ belongs to the set
$$
\{2k:\ k=0,\pm1,\ldots, \pm 2\beta \},
$$ 
because $\gamma_0-\gamma'_{2\beta+1}\in \Gamma-\Gamma$. However,  this leads to a contradiction too. On the one hand,  in view of \eqref{contra1} and  
\eqref{contra2} 
\begin{align*}
\gamma_0+\la_0&=(a_k,k)\ \\
\gamma'_{2\beta+1}+\la'_{2\beta+1}&=(2\beta+1,k-1)\,
\end{align*}
and so  we have $(\gamma_0-\gamma_{2\beta+1}) + (\la_0-\la_{2\beta+1})=(a_k-(2\beta+1),1)$, which gives
$$
(s,\mu):= \la_0-\la_{2\beta+1}'=(a_k-(2\beta+1),1-n).
$$
On the other hand, by Theorem \ref{geom cond} we obtain that $M^{-1}(s,\mu)\in \cZ (V_g g)$. However, 
$$
|V_gg(M^{-1}(s,\mu))|=|V_gg(s/2,\mu/(2\beta+1))|=|\widehat{\chi}_{\W\cap(s/2+\W)}(\mu/(2\beta+1))|\neq 0
$$
because 
$$
\W\cap(s/2+\W)=\left[ \frac{a_k}{2}, \frac{1}{2} \right),
$$
which has measure less than $\frac12$, and $|\mu/(2\beta+1)|<2$ (again we have used Lemma  \ref{easy}).
Thus, this completes the proof of the claim. 

\medskip

Therefore, we conclude that
\begin{equation} \label{1812}
\Gamma+\La_M=\bigcup_{k\in \Z} \{k\} \times  (\Z+a_k),  \
\end{equation}
where $a_k$ are real numbers in $[0,1)$ for $k\neq 0$, and $a_0=0$. Our next step will be to prove that
\begin{equation} \label{marcha eslava}
\La_M=\bigcup_{k\in K} \{k\} \times  (L_k+a_k)
\end{equation}
where $K$ is a tiling complement of $\{0,2\beta+1\}$ in $\Z$, and $L_k\subseteq \Z$ whose structure will be studied later. With this aim, it is enough to prove that for every $k,n\in\Z$, if
\begin{align*}
(k,a_k+n)=\gamma_0+\la_0, \peso{and}
(k,a_k+n+1)=\gamma_1+\la_1,
\end{align*}
then $(\gamma_0)_1=(\gamma_1)_1$. Suppose that it is not the case, hence $ (\gamma_0-\gamma_1)_1=\pm(2\beta+1)$. In consequence, we obtain that
$$
\la_0-\la_1=(\mp(2\beta+1), m),
$$
where $|m|< 4\beta+2$. As in the proof of the claim, this leads to a contradiction because $M^{-1}(\la_j-\la_k)$ does not belongs to the zero set of $V_gg$.
This proves that
\begin{align*}
\{0,2\beta+1\}+K=\Z &\Longleftrightarrow ([0,1)\cup[2\beta+1,2\beta+1))+K=\R\\
&\Longleftrightarrow ([0,1/2)\cup[\beta+1/2,\beta+1))+\frac{K}{2}=\R.
\end{align*}
This implies that $(\W+k/2)\cap(\W\cap k'/2)=\varnothing$ for any pair of different elements $k,k'\in K$. Since $\{g_{\la}\}_{\la \in \La}$ is a Gabor orthonormal basis of $L^2(\R)$, we get that for each $k$, the set $\frac{L_k+a_k}{2\beta+1}$ is a spectrum for $\W+k/2$. This is equivalent to saying that for each $k\in K$, the sets $\frac{L_k}{2\beta+1}$ are spectra for $\W$. So, to conclude the proof, it is enough to prove that 
$$
A:=\W\peso{and }  B:=\bigcup_{k=0}^{2\beta} \left[\frac{2k}{2\beta+1}, \frac{2k+1}{2\beta+1}\right )
$$
are (tight) spectral pairs (see Definition \ref{Tannhauser}). Indeed, if these two sets are spectral pairs, by Theorem \ref{KOO}, the set $B$ tiles $\R$ with 
$\frac{L_k}{2\beta+1}$. But, this is equivalent to saying that 
$$
\bigcup_{k=0}^{2\beta} \left[{2k}, {2k+1}\right )
$$
tiles the real line with $L_k$, or equivalently $L_k+\{2n:\ n=0,\pm1,\ldots,\pm2\beta\}=\Z$. Since $|A|=|B|=1$, it is enough to prove that $A$ and $B$ are spectral pairs, which by definition means that:
\begin{itemize}
\item[a.)] $(B^{\circ}-B^{\circ}) \cap \cZ(\widehat{\chi}_{A})=\varnothing$;
\item[b.)] $(A^{\circ}-A^{\circ}) \cap \cZ(\widehat{\chi}_{B})=\varnothing$
\end{itemize}

\medskip

As in Lemma \ref{dif2}, we can prove that in the interval $(-2,2)$, the unique zeros of $\widehat{\chi}_A$ are those of the form 
$\w=\frac{k}{2\beta+1}$  where $k$ is an odd integer. Therefore, we get (a). On the other hand, 
$A^{\circ}-A^{\circ}=(-1/2,1/2)\cup (\beta,\beta+1)\cup (-\beta-1,-\beta)$, and straightforward computations show that
\begin{align*}
|\widehat{\chi}_B(\w)|
&= 
\left|\frac{\sin \frac{\pi \w}{2\beta+1}}{\pi \w}\,\cdot \,\sin 2\pi \w\,\cdot \,\frac{1}{\sin \frac{2\pi \w}{2\beta+1}}\right|.
\end{align*}

Note that none of the three sines vanish at $(-1/2,1/2)\setminus\{0\}$, and clearly zero is not a problem because $\widehat{\chi}_B(0)=|B|=1$ . The other points to take into account are $\pm(\beta+1/2)$. At these points, the last two sines cancel each other and the other part of the expression does not vanish. In consequence, (b) also holds, and the sets $A$ and $B$ are spectral pairs. 


%
%
%

\appendix
{	
\section{Appendix}
\subsection*{Description of the zero set of $V_{\chi_\W}\chi_\W$ when $\W$ tiles $\R$} \label{Zero set}

Recall the set
$$
\W=[0,\alpha] \cup [\beta+\alpha, \beta+1].
$$
 Throughout this section we completely describe the zero set of the 
Short time Fourier transform of $g=\chi_\W$ in each of the following cases:
 \begin{itemize}
\item[(i)]  $\displaystyle 0<\alpha < 1/2 $ and  $\beta \in \N$; 
\item[(ii)] $\alpha =1/2$ and $\displaystyle \beta \in \frac{1}{2}\N$. 
\end{itemize}

The results of this section have not been necessary to obtain the results of the previous sections. However, the detailed description that follows may be useful because it clearly encodes the orthogonality of the time-frequency translates.

\medskip

Recall that the zero set of $V_gg$, is given by 
$$
\cZ(V_gg)=\{(t,\nu): V_gg(t,\nu)=0\}.
$$
By the symmetries of this set, due to Lemma \ref{symmetries}, it is enough to study the subset 
$$
\cZ^+(V_gg)=\{(t,\nu):\ (t,\nu)\in \cZ(V_gg)\ \ t,\nu\geq 0\}.
$$

\medskip

 As we mentioned before, $V_gg(t,\nu)=\widehat{\chi}_{\W\cap(\W+t)}(\nu)$. If $t\in[0,\alpha)$ then
\begin{equation}\label{eq two intervals}
\W\cap(\W+t)=[t,\alpha]\cup[\beta+\alpha+t,\beta +1],
\end{equation}
while if $t>\alpha$ the set $\W\cap(\W+t)$ is a single interval. This gives a different structure of the zeros, depending on the value of $t$. 
Hence, we will divide the study of (i) and (ii) in two cases. Let  $\cZ^+(V_gg)=\cZ^+_1(V_gg)\cup \cZ^+_2(V_gg)$,  where
\begin{align*}
\cZ^+_1(V_gg)&=\{(t,\nu):\ (t,\nu)\in \cZ(V_gg)\ t\in[0,\alpha)\ \ \nu\geq 0\};\\
\cZ^+_2(V_gg)&=\{(t,\nu):\ (t,\nu)\in \cZ(V_gg)\ t\in[\alpha,+\infty)\ \ \nu\geq 0\}.\\
\end{align*}

\subsection*{The case  $\displaystyle 0<\alpha < 1/2$ and $\beta \in \N$}

%


As a direct consequence of Lemma \ref{easy} we get that

$$
\cZ^+_2(V_gg)=\begin{cases}
   (t,\nu)                                       &\mbox{$t \in [1-\alpha, \beta]\cup [\beta+1,\infty),  \nu >0$ }\\
   (t,{k}/{(1-\alpha-t)})                   &\mbox{$t\in[\alpha,1-\alpha), k\in \N $}\\ 
   (t,{k}/{(t-\beta)})                        &\mbox{$t\in[\beta, \beta+\alpha),  k\in \N$}\\
   (t,k/{\alpha})                             &\mbox{$t\in[\beta+\alpha, \beta+1-\alpha), k\in \N $}\\
   (t,k /(\beta+1-t))                        &\mbox{$t\in[\beta+1-\alpha, \beta+1), k\in \N$}.\\
\end{cases} 
$$

\medskip

On the other hand, the following result describes $\cZ^+_1(V_gg)$.

\begin{pro} \label{zeroes1}
Let $g=\chi_\W$. Then
 $$
\cZ^+_1(V_gg)=\begin{cases}
  \displaystyle \Big(\alpha-\frac{k}{\nu}\,,\,\nu\Big)     
  &\mbox{$n,k\in \N$, $\displaystyle \nu=\frac{n}{1-2\alpha}$,  and $\displaystyle \frac{k}{\nu}\in (0,\, \alpha]$}\\
  \ & \ \\
  \displaystyle \Big(\frac{k}{n}- \beta,n\Big)                                 
  &\mbox{$ k,n\in\N$ and $\displaystyle\frac{k}{n}\in [\beta,\beta+ \alpha)$}.\\
\end{cases} 
$$ 
 If in addition, there exists $r\in \Q$ of the form $\displaystyle \frac{2k_0+1}{2n_0}$, so that $\beta= r-(2r+1)\alpha$, then the set 
$$
\Big\{ \Big(t, \frac{n}{1-2\alpha} \Big) : t\in [0, \alpha),\,  n\in \N \, \mbox{such that  $\displaystyle r=\frac{2k+1}{2n}$ for some $k\in\N$ }\Big\}
$$
should be added to the above zero set. 
\end{pro}

\bdem[Proof of Proposition \ref{zeroes1}]

First of all, note that we have that $V_gg(t, 0)\neq 0$ if $t\in[0,\alpha)$. Indeed, this follows directly by \eqref{Vg} because
$\W\cap(\W+t)$ has non-empty interior. So, from now on we will assume that $\nu>0$. As we observed in \eqref{eq two intervals},  
$$
\W\cap (\W+t)=[t,\alpha]\cup[\beta+\alpha+t,\beta +1]
$$ 
Then, a direct computation shows that $V_gg(t,\nu)= 0$ if 
$$
e^{-2\pi i\nu t}-e^{-2\pi i\nu \alpha} +e^{-2\pi i\nu (\beta+\alpha+t)} -e^{-2\pi i\nu (\beta+1)} =0
$$
Since, for given $z_i,$ so that $|z_i|=1$ for every $i=1,2,3,4$, all solutions to $z_1+z_2+z_3+z_4=0$ are given by any
pairs of opposite numbers and so we have the following cases: 
\begin{itemize}
\item[(i)] $e^{-2\pi i\nu t}=e^{-2\pi i\nu \alpha}$ and $e^{-2\pi i\nu (\beta+\alpha+t)} = e^{-2\pi i\nu (\beta+1)}$;
\item[(ii)] $e^{-2\pi i\nu t}=e^{-2\pi i\nu (\beta+1)}$  and $e^{-2\pi i\nu \alpha} =e^{-2\pi i\nu (\beta+\alpha+t)} $;
\item[(iii)] $-e^{-2\pi i\nu t}= e^{-2\pi i\nu (\beta+\alpha+t)}$ and $-e^{-2\pi i\nu \alpha} = e^{-2\pi i\nu (\beta+1)} $.
\end{itemize}

\begin{description}

\item[Case (i):] In this case we have that
$$
\nu (\alpha-t)= k \quad \text{ and } \quad \nu (1-2\alpha)= n,
$$
where $k,n\in \N$. So  
$$
\nu= \frac{n}{1-2\alpha}  \quad \text{and} \quad t=\alpha + (2\alpha-1)\frac{k}{n}.
$$
Since $t\in [0,\alpha)$ we have that $k/n \in (0, \alpha/(1-2\alpha)]$. 

\item[Case (ii):] In this case
$$
\nu (\beta+t)=k \quad \text{ and } \quad \nu (t+\beta+1) =m,
$$
where $k,m\in \N$.  This yields
$$
\nu (\beta+t)=k \quad \text{ and } \quad \nu=n \in\N.
$$
and so $\displaystyle t=\frac{k}{n}-\beta$. Since $t\in [0,\alpha)$ we have that $k/n \in [\beta,\beta+\alpha)$.

\item[Case (iii):] Finally, in this case we get that 
$$
2\nu(\beta+\alpha)=2k+1 \quad \text{ and } \quad 2\nu (1+\beta-\alpha)=2m+1
$$
where $k,m \in \N$. Hence
$$
2\nu(\beta+\alpha)=2k+1 \quad \text{ and } \quad \nu (1-2\alpha)=n
$$
with $n\in \N$. So, in this case has solutions only if 
$$
\frac{\beta+\alpha}{1-2\alpha}= \frac{2k+1}{2n}=r\in \Q,
$$
hence the system can be solved if and only if $\beta = r-(2r+1) \alpha$. So, we get the additional case.
\end{description}
\edem

\subsection*{Example}\label{Exa1}

Take $\alpha=1/ \sqrt 15$  and $\beta=2$, and let $g=\chi_\W$.  Note that we choose $\alpha, \beta$ 
to be rationally independent, because in this case the zero set of $V_gg$ is simpler. Indeed, in this case the set $\cZ^+_1(V_gg)$ 
is described by a union of two sets. Otherwise, we may have to consider one more case, as it is shown in Proposition \ref{zeroes1}.

\medskip

As we have seen the set $\cZ^+(V_gg)$ has different structure depending on the value of $t$.  
When $t\in[1/ \sqrt 15, \infty)$, the set $\cZ^+_2(V_gg)$ is described by the following picture.

\begin{figure}[H]
           \centering
           \includegraphics[height=5cm]{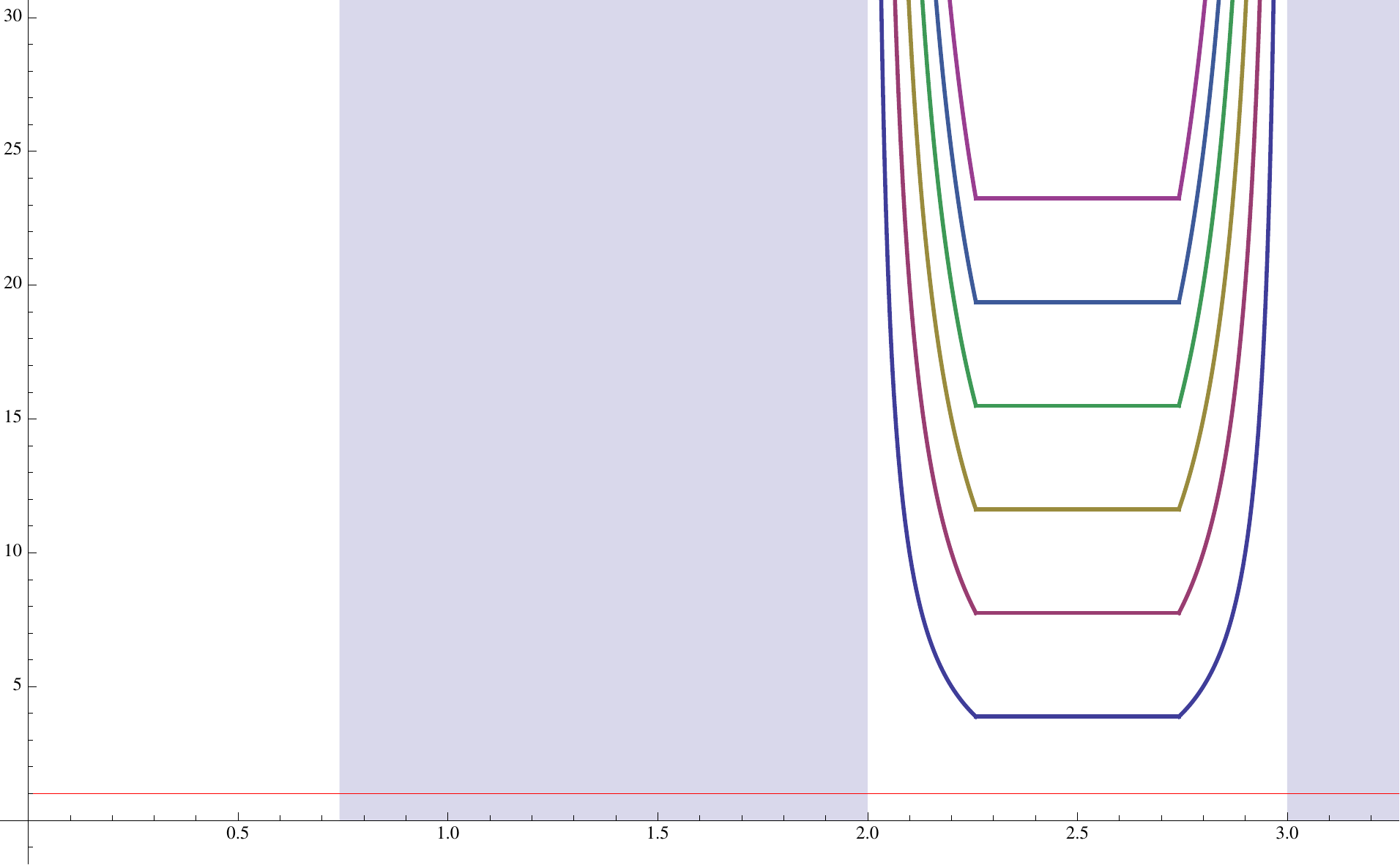}
           \par\vspace{0cm}
           \caption{The set $\cZ^+_2(V_gg)$.} 
           \label{z_2}
\end{figure}

When $t\in[0,1/ \sqrt 15)$ we get the set $\cZ^+_1(V_gg)$, which is described by Proposition \ref{zeroes1}, as union of two sets.  
The first one corresponds to the set
$$
 \Big(\frac{1}{\sqrt 15}-\frac{k}{\nu}\,,\,\nu\Big)     
$$
where $n,k\in \N$, $\displaystyle \nu=\frac{n}{1-2/ \sqrt 15}$,  and $\displaystyle \frac{k}{\nu}\in (0,\, 1/ \sqrt 15]$.

\begin{figure}[H]
           \centering
           \includegraphics[height=4.5cm]{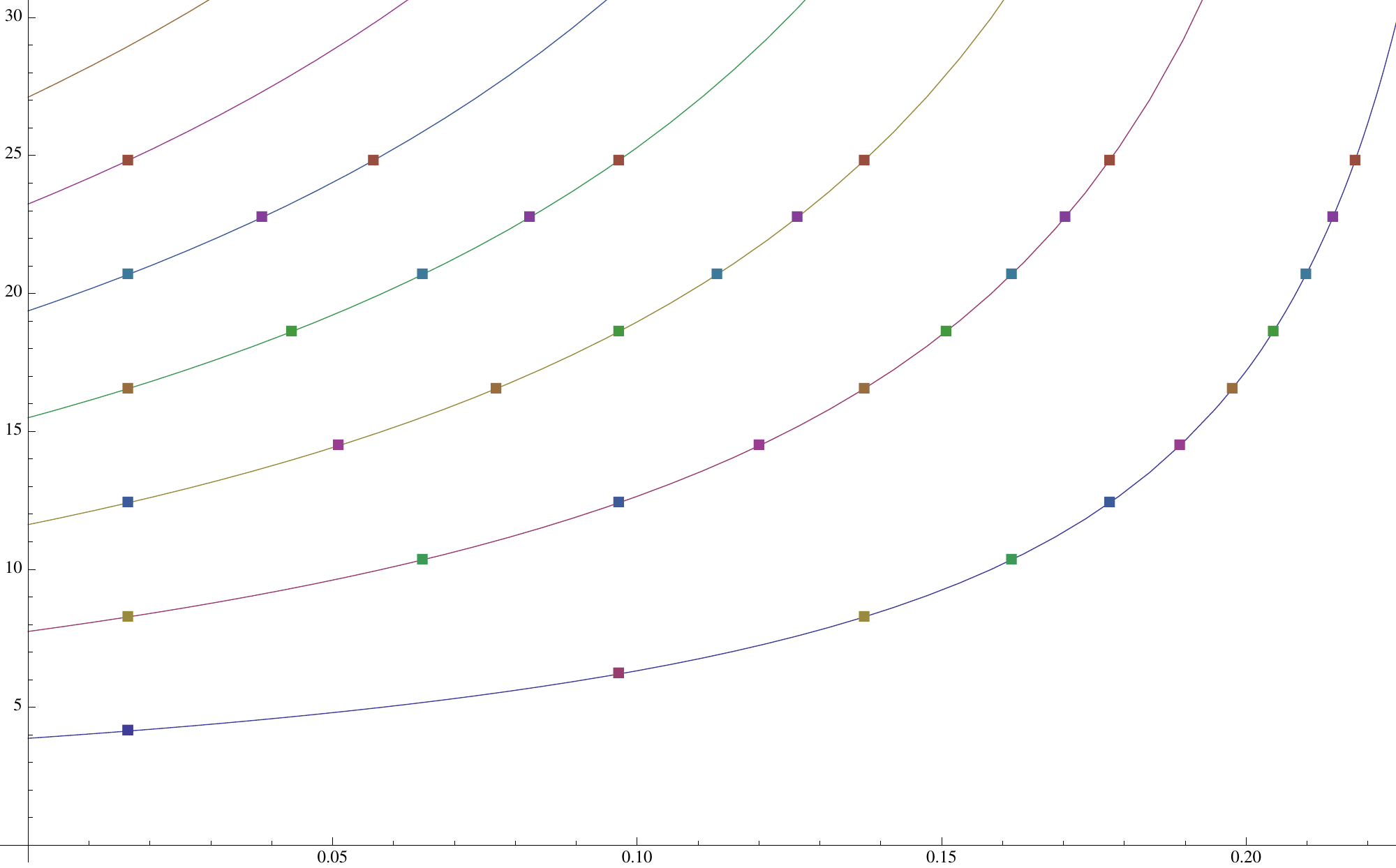}
           \par\vspace{0cm}
           \caption{The squares describe the first subset of $\cZ^+_1(V_gg)$ in Proposition \ref{zeroes1}.}
           \label{z_11}
\end{figure}

The second one corresponds to the set 
$$
\Big(\frac{k}{n}- \beta,n\Big),                            
$$
 where $ k,n\in\N$ and $\displaystyle\frac{k}{n}\in [\beta,\beta+ \alpha)$. 
 
\begin{figure}[H]
           \centering
           \includegraphics[height=5cm]{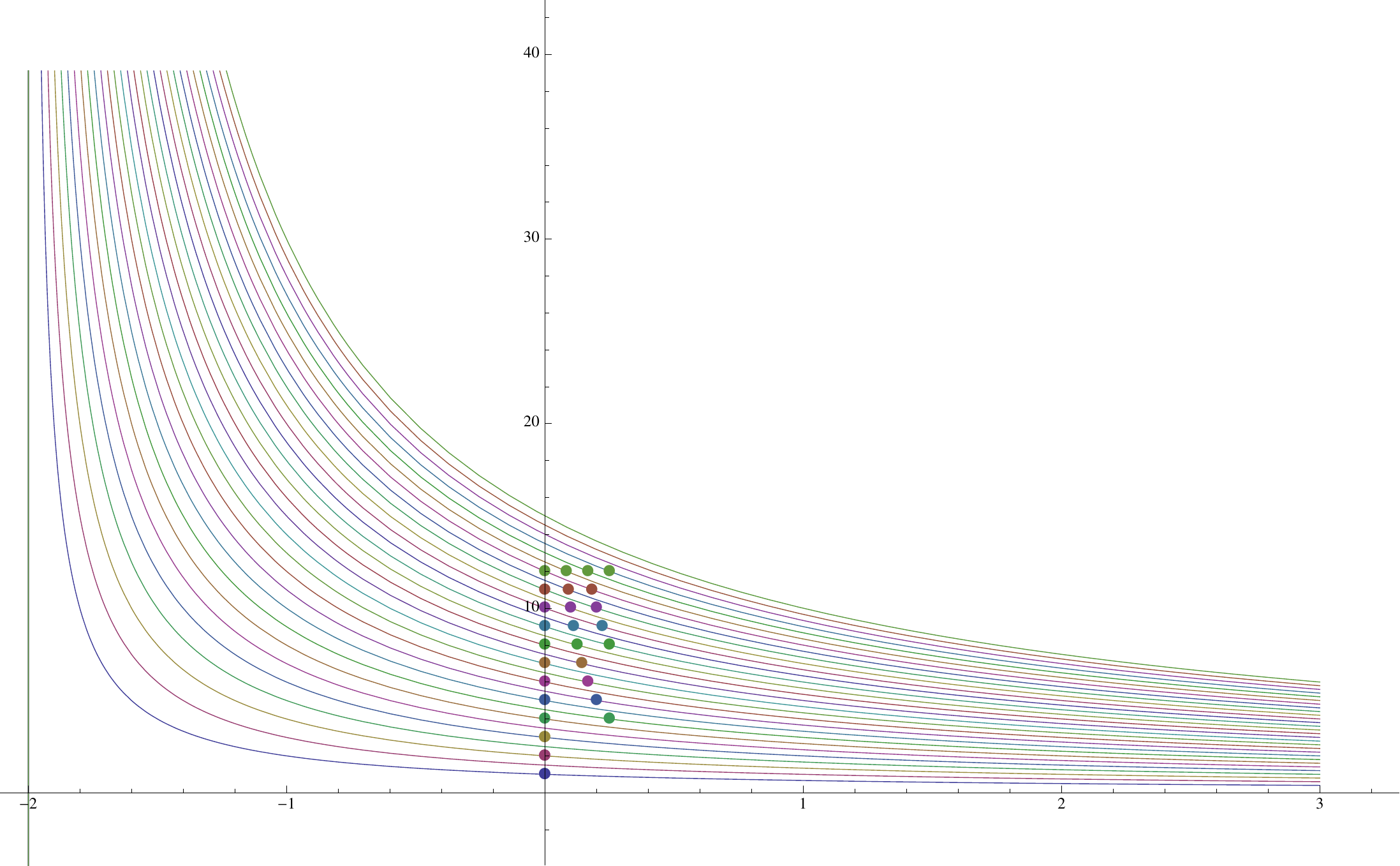}
           \par\vspace{0cm}
           \caption{The circles describe the second subset of $\cZ^+_1(V_gg)$  in Proposition \ref{zeroes1}.}
           \label{z_12}
\end{figure}

Figure \ref{z_all} describes completely the set $\cZ^+(V_gg)$. 

\begin{figure}[H]
           \centering
           \includegraphics[height=5cm]{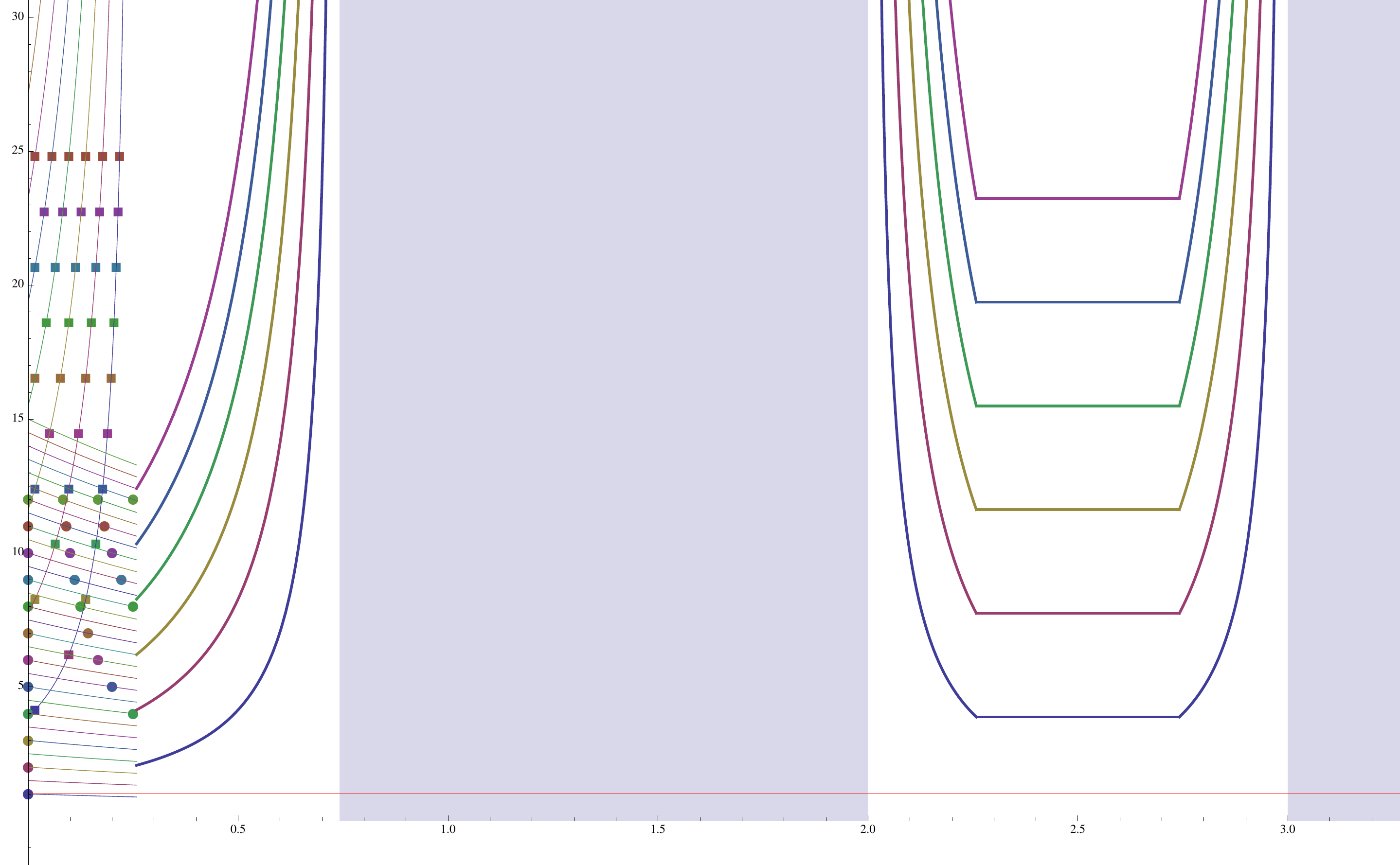}
           \par\vspace{0cm}
           \caption{The set $\cZ^+(V_gg)$ corresponding to $\alpha=1/ \sqrt 15$ and $\beta=2$.}
           \label{z_all}
\end{figure}

\subsection*{The case $\alpha =1/2$ and $\displaystyle  \beta \in \frac{1}{2}\N$}

%

As in the case $\alpha<1$, the description of $\cZ^+_2(V_gg)$ is much simpler. Indeed, as a direct consequence of Lemma  \ref{easy} we get that

$$
\cZ^+_2(V_gg)=\begin{cases}
   (t,\nu)                                       &\mbox{$t \in [1/2, \beta]\cup [\beta+1,\infty),  \nu >0$ }\\
   (t,{k}/{(t-\beta)})                        &\mbox{$t\in (\beta, \beta+1/2),  k\in \N$}\\
   (t,k /(\beta+1-t))                        &\mbox{$t\in [\beta+1/2, \beta+1), k\in \N$}.\\
\end{cases}
$$

\medskip

On the other hand, the description of $\cZ^+_1(V_gg)$ is given in the following proposition.

\begin{pro} \label{halfzeroes1}
Let $\W=[0,1/2) \cup [1/2+\beta,1+ \beta)$, and $g=\chi_\W$. Then if $t \in [0,1/2]$
$$
\cZ^+_1(V_gg)=\begin{cases}
 \displaystyle \Big(t, \frac{k}{1/2-t}\Big)     
  &\mbox{$t\in [0,1/2]$, and $k\in \N$}\\
  \ & \ \\
  \displaystyle \Big(\frac{k}{\nu}-\beta\,,\,\nu\Big)     
  &\mbox{$k\in \N$, $\displaystyle \nu=\frac{n}{2\beta+1}$,  and 
  $\displaystyle \frac{k}{n}\in \Big[\frac{\beta}{2\beta+1},\, \frac{1}{2}\Big]$}\\
  \ & \ \\
  \displaystyle \Big(t, \frac{2k}{2\beta +1}\Big)  &\mbox{$ k\in\N$ and $t\in [0,1/2]$}.\\
\end{cases} 
$$
\end{pro}

\bdem
Recall that in \eqref{intersectionshalf} we have
\begin{equation} \label{intersectionsa}
\W\cap(\W+t)=\begin{cases}
[t,1/2]\cup[\beta+1/2+t,\beta +1]&\mbox{if $t\in[0,1/2)$}\\
[\beta+1/2, t+1/2]&\mbox{if $t\in[\beta,\beta+1/2)$}\\
[t,\beta+1]&\mbox{if $t\in[\beta+1/2,\beta+1)$}.
\end{cases}
\end{equation}

Then, a direct computation shows that $V_gg(t,\nu)= 0$  and  $t \in [0,1/2]$  we have that the following system should be satisfied:
$$
e^{-2 \pi i \nu t}- e^{-\pi i \nu }+ e^{-2\pi i \nu(\beta+ 1/2+ t)} - e^{-2\pi i \nu (\beta +1)}= 0 
$$
Since, for given $z_i,$ so that $|z_i|=1$ for every $i=1,2,3,4$, all solutions to $z_1+z_2+z_3+z_4=0$ are given by any
pairs of opposite numbers and so we have the following cases: 
\begin{itemize}
\item[(i)] $e^{-2\pi i\nu t}=e^{- \pi i\nu}$ and $e^{-2\pi i\nu (\beta+1/2+t)} = e^{-2\pi i\nu (\beta+1)}$;
\item[(ii)] $e^{-2\pi i\nu t}=e^{-2\pi i\nu (\beta+1)}$  and $e^{-\pi i\nu} =e^{-2\pi i\nu (\beta+1/2+t)} $;
\item[(iii)] $-e^{-2\pi i\nu t}= e^{-2\pi i\nu (\beta+1/2+t)}$ and $-e^{-\pi i\nu} = e^{-2\pi i\nu (\beta+1)} $.
\end{itemize}

\begin{description}

\item[Case (i):] In this case we have that for $t\in [0,1/2]$ 
$$
\nu (1/2-t)= k, 
$$
where $k \in \N$. 

\item[Case (ii):] In this case
$$
\nu (t-\beta-1)\in \N \quad \text{ and } \quad \nu (t+\beta) \in \N,
$$
where $k,m\in \N$.  This yields
$$
\nu (\beta+t)=k \quad \text{ and } \quad \nu(2\beta+1)= n.
$$
So $\displaystyle t=\frac{k}{\nu}-\beta$ and $\displaystyle\nu =\frac{n}{2\beta+1}$.  
Since $t\in [0,1/2)$ we have that $\displaystyle \frac{k}{n} \in \Big[\frac{\beta}{2\beta+1},\,  \frac{1}{2}\Big]$.

\item[Case (iii):] Finally,  we get that 
$$
2\nu(\beta+1/2)=2k+1 
$$
where $k \in \N\cup \{0\}$.
\end{description}
\edem

\subsection*{Example} \label{Exa2} 

Let $\alpha=1/2$ and $\beta=2$.  The next picture corresponds to $\cZ^+_2(V_gg)$
\begin{figure}[H]
           \centering
           \includegraphics[height=5cm]{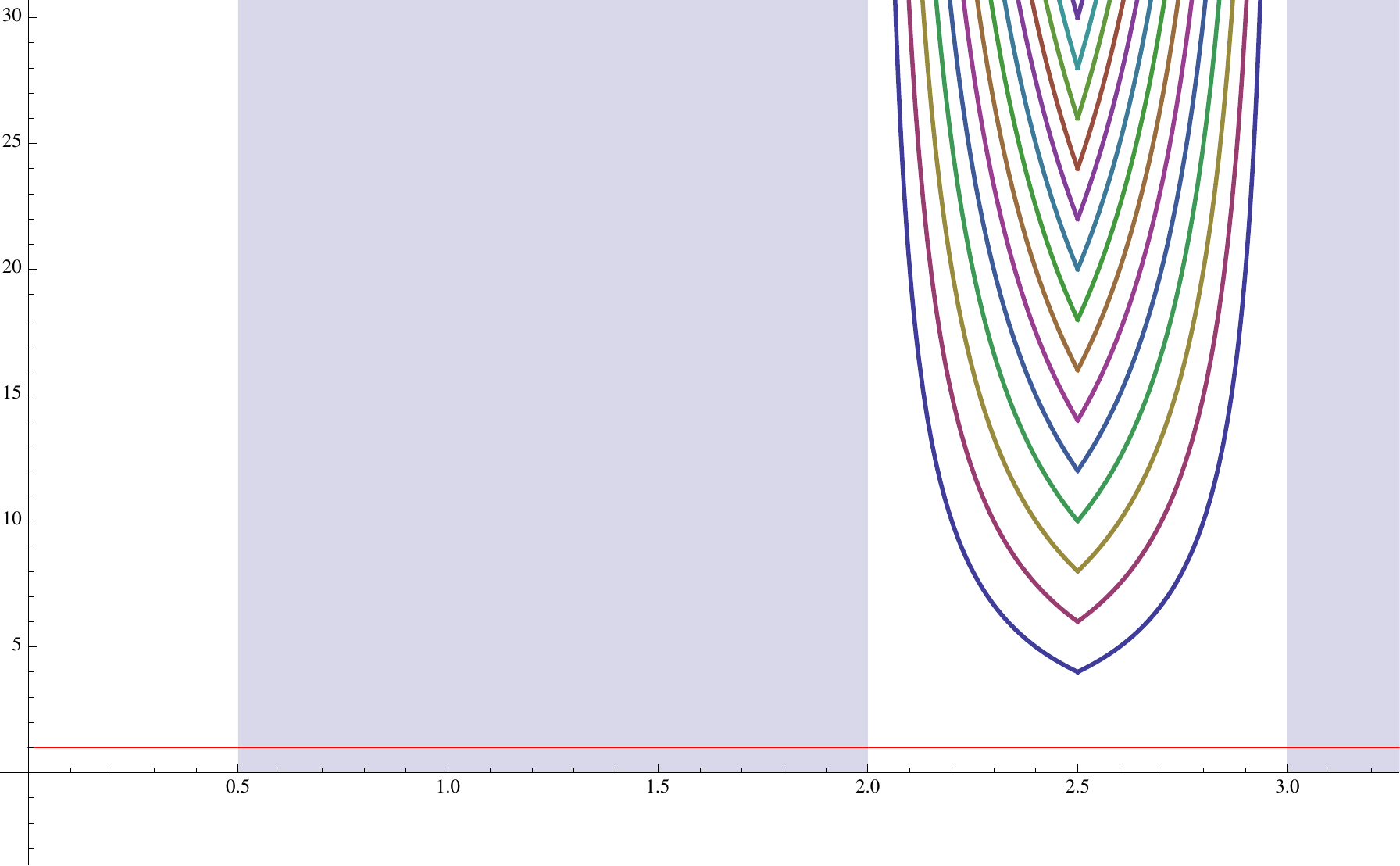}
           \par\vspace{0cm}
           \caption{The set $\cZ^+_2(V_gg)$.} 
           \label{ahalf_z2}
\end{figure}

The set $\cZ^+_1(V_gg)$, is given as union of three sets as calculated by Proposition \ref{halfzeroes1}.  Any of this set is described by figures
\ref{ahalf_z11}, \ref{ahalf_z12} and \ref{ahalf_z13}.

\begin{figure}[H]
           \centering
           \includegraphics[height=5cm]{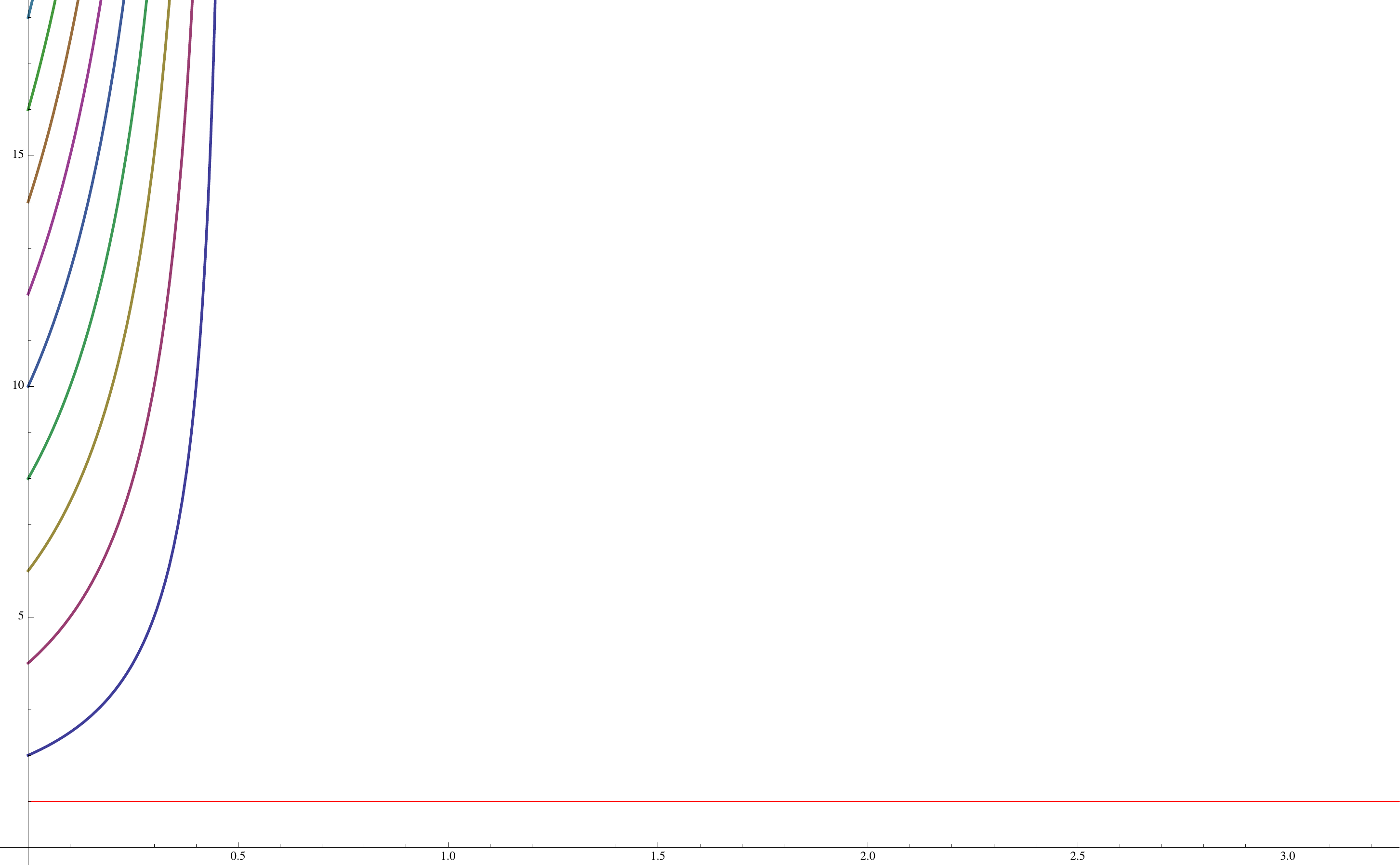}
           \par\vspace{0cm}
           \caption{The hyperbolas describe the first subset of $\cZ^+_1(V_gg)$ in Proposition \ref{halfzeroes1}.}
           \label{ahalf_z11}
\end{figure}

\medskip

\begin{figure}[H]
           \centering
           \includegraphics[height=5cm]{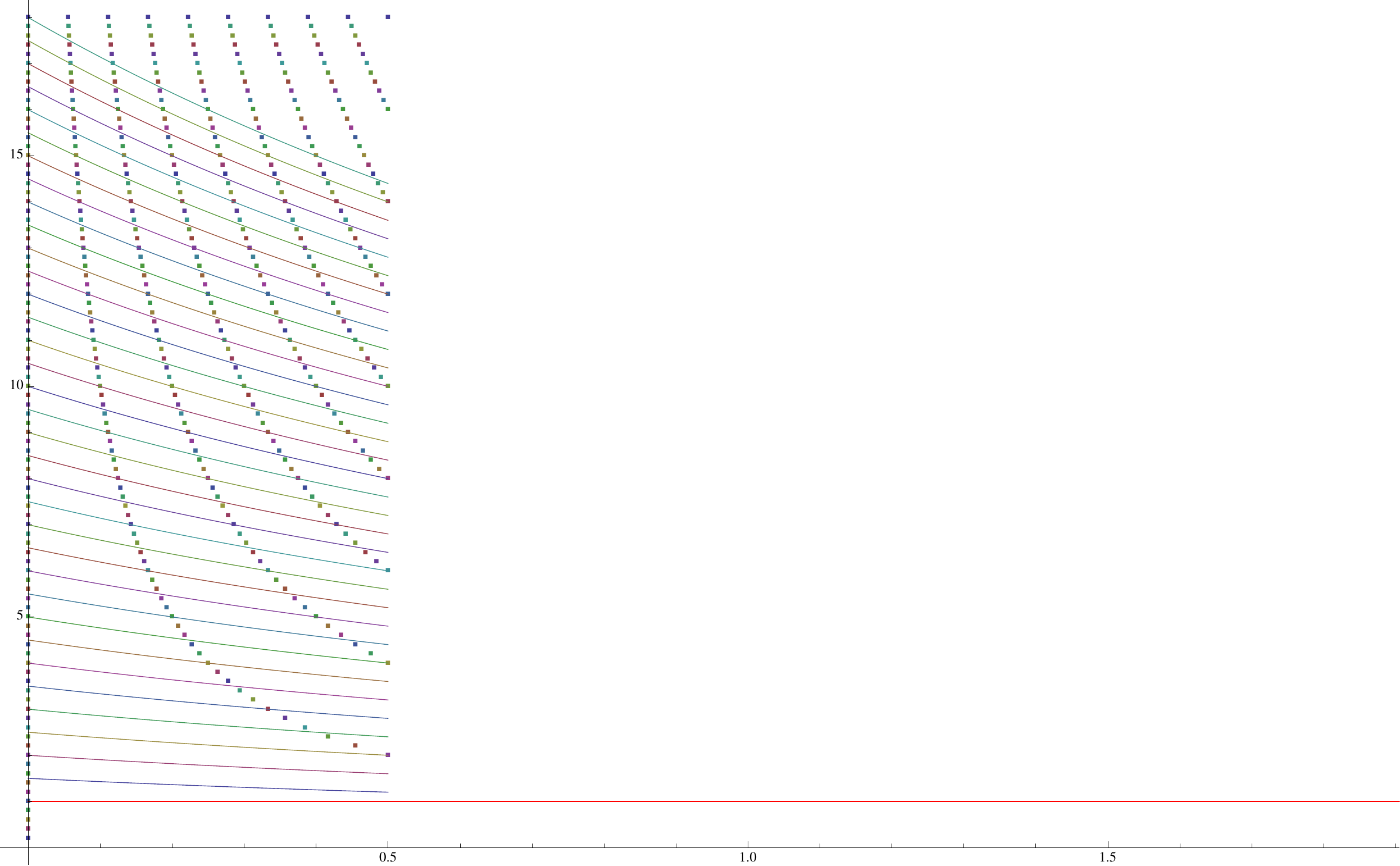}
           \par\vspace{0cm}
           \caption{The dots describe the second subset of $\cZ^+_1(V_gg)$  in Proposition \ref{halfzeroes1}.}
           \label{ahalf_z12}
\end{figure}

\medskip

\begin{figure}[H]
           \centering
           \includegraphics[height=5cm]{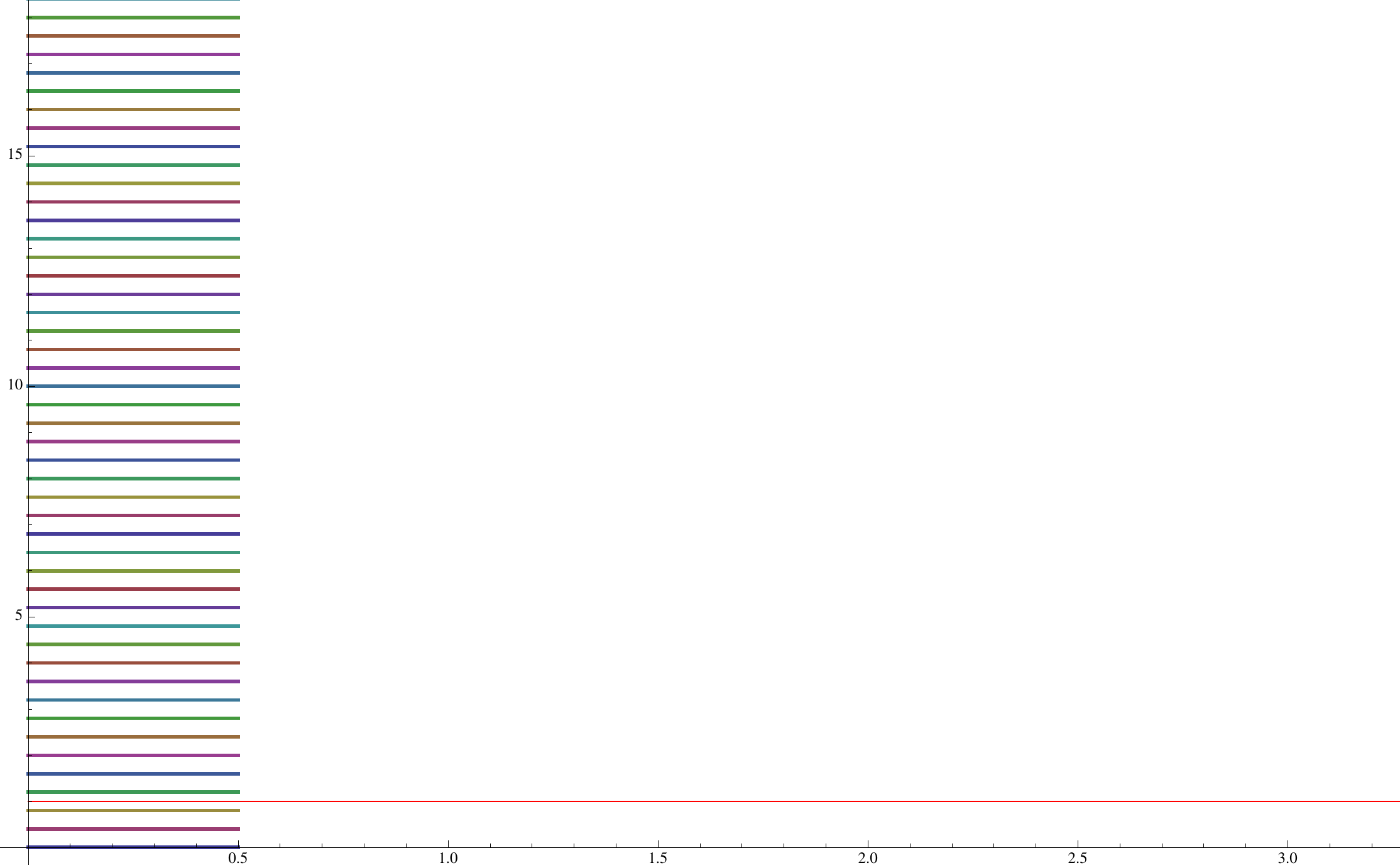}
           \par\vspace{0cm}
           \caption{The lines describe the third subset of  $\cZ^+_1(V_gg)$  in Proposition \ref{halfzeroes1}.}
           \label{ahalf_z13}
\end{figure}

\medskip

Finally the set $\cZ^+(V_gg)$ is given in Picture \ref{Zeroes_1}.

\medskip
\begin{figure}[H]
           \centering
           \includegraphics[height=5cm]{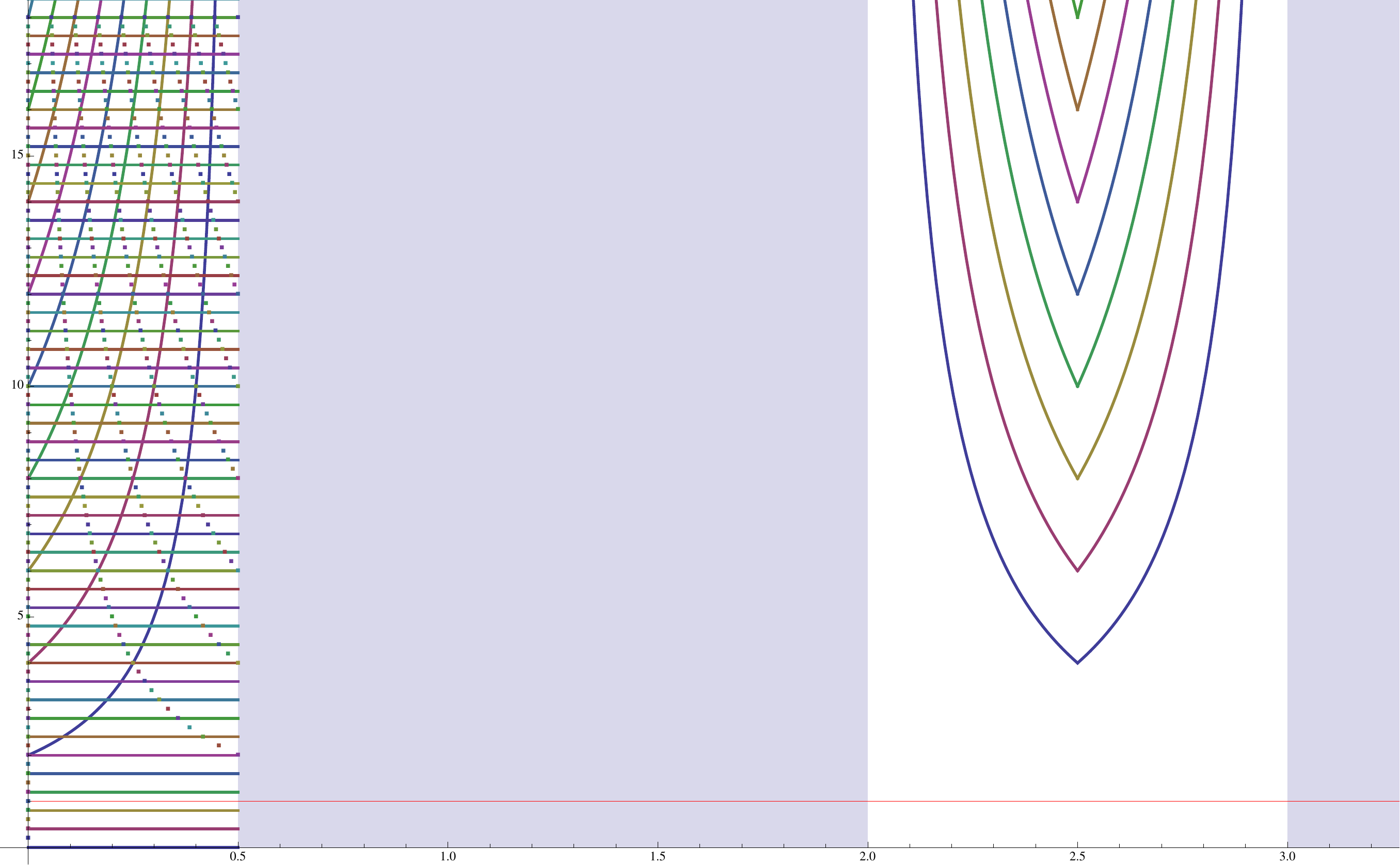}
           \par\vspace{0cm}
           \caption{The set $\cZ^+(V_gg)$  corresponding to the case $\alpha=1/2$ and $\beta=2$.}
           \label{Zeroes_1}
\end{figure}
 
}

\bigskip

\textsc{Elona Agora,}
Instituto Argentino de Matemática \lq \lq Alberto P. Calderón" (IAM-CONICET), Buenos Aires, Argentina \\
\hfill\textit{E-mail address: elona.agora@gmail.com} 

\bigskip

\textsc{Jorge Antezana,}
Departamento de Matemática, Universidad Nacional de La Plata and, 
 Instituto Argentino de Matemática  \lq \lq Alberto P. Calderón" (IAM-CONICET), Buenos Aires, Argentina \\
 \hfill  \textit{E-mail address: antezana@mate.unlp.edu.ar}

\bigskip

\textsc{Mihail N. Kolountzakis,}
Department of Mathematics and Applied Mathematics
University of Crete, Heraklion, Greece \\
 \hfill  \textit{E-mail address: kolount@gmail.com}

\end{document}